\documentclass[12pt]{amsart}
\usepackage[dvipsnames]{xcolor}
\usepackage{latexsym}
\usepackage{a4wide}
\usepackage{amscd}
\usepackage{graphics}
\usepackage[pdftex]{graphicx}
\usepackage{amsmath}
\usepackage{amssymb}
\usepackage{esint}
\usepackage{mathrsfs}
\usepackage{mathabx} 
\usepackage{hyperref}
\usepackage{amsthm}
\usepackage{bbm}
\usepackage{color}
\usepackage{accents}
\usepackage{enumerate}
\usepackage{enumitem}
\input xy
\xyoption{all}

\usepackage{tikz-cd}
\usepackage{tikz,float,caption}
\usetikzlibrary{patterns} 
\usetikzlibrary{arrows.meta,calc,decorations.markings,patterns,cd,patterns.meta}

\usetikzlibrary{decorations.pathreplacing}

\newcommand{\Z}{\mathbb{Z}}

\newcommand{\R}{\mathbb{R}}
\newcommand{\C}{\mathbb{C}}

\newcommand{\F}{\mathbb{F}}
\newcommand{\D}{\mathbb{D}}

\newcommand{\eps}{\epsilon}
\newcommand{\CF}{\mathrm{CF}}
\newcommand{\HF}{\mathrm{HF}}
\renewcommand{\H}{\mathrm{H}}
\renewcommand{\d}{\mathrm{d}}
\newcommand{\varphitilde}{\widetilde{\varphi}}

\newcommand{\spec}{\mathcal{S}}             
\newcommand{\specm}{\mathcal{S}_{\mathrm{mean}}}            

\newcommand{\A}{\mathbb{A}}
\newcommand{\wind}{\mathrm{wind}}  	
\newcommand{\Fix}{\mathrm{Fix}}		

\newcommand{\Ham}{\mathrm{Ham}}

\newcommand{\Hof}{\mathrm{Hof}}
\newcommand{\red}[1]{\textcolor{red}{#1}}

\newtheorem{theorem}{Theorem}[section]
\newtheorem{corollary}[theorem]{Corollary}
\newtheorem{claim}[theorem]{Claim}
\newtheorem{definition}[theorem]{Definition}
\newtheorem{example}[theorem]{Example}
\newtheorem{question}[theorem]{Question}
\newtheorem{remark}[theorem]{Remark}
\newtheorem{lemma}[theorem]{Lemma}
\newtheorem{proposition}[theorem]{Proposition}

\author[H. Alizadeh]{Habib Alizadeh}
\address{Habib Alizadeh, The Institute of Geometry and Physics, University of Science and Technology of China, 96 Jinzhai Road, Baohe District, Hefei, Anhui, 230000, China}
\email{\texttt{habib.alizadeh.math@gmail.com}}

\author[D. Cristofaro-Gardiner]{Dan Cristofaro-Gardiner}
\address{Dan Cristofaro-Gardiner, University of Maryland, 4176 Campus Drive - William E. Kirwan Hall
College Park, MD 20742}
\email{\texttt{dcristof@umd.edu}}

\author[A. Pirnapasov]{Abror Pirnapasov}
\address{Abror Pirnapasov, New Uzbekistan University, Movarounnahr street 1, 100007 Tashkent, Uzbekistan}
\email{\texttt{a.pirnapasov@newuu.uz}}

\author[E. Shelukhin]{Egor Shelukhin}
\address{Egor Shelukhin, Department of Mathematics, University of Montreal, Canada}
\email{\texttt{egorshel@gmail.com}}

\title{The topology of the mean action spectrum}
\begin{document}

\begin{abstract}
We prove that if an area-preserving diffeomorphism of the closed two-disk is neither a pseudo-rotation nor a periodic map, then the closure of its mean action spectrum contains an interval of positive length. This partially answers a question of Hutchings. We do not require any additional assumptions on the area-preserving map, such as non-degeneracy or specific behavior near the boundary.
\end{abstract}

\maketitle
\tableofcontents


\section{Introduction}
Associated to a Hamiltonian diffeomorphism of an aspherical surface $(\Sigma, \sigma)$ there is a natural closed nowhere dense subset $\spec(\phi) \subset \mathbb R,$ called the action spectrum. If $\Sigma$ has boundary with a distinguished boundary component, then $\sigma$ is exact. Given a choice of primitive $\lambda$ of $\sigma,$ the spectrum is defined by choosing a generating time-dependent Hamiltonian $H$ vanishing on a fixed boundary component, and looking at the critical values of the action functional \[ \mathcal A_H(z) = \int_0^1 H(t, z(t))\,dt - \int_z \lambda\] where $z: \R/\Z \to \Sigma$ is a contractible loop in $\Sigma.$ It is not hard to see that this function only depends on a choice of isotopy $(\phi^t)$ from the identity to $\phi.$

On the other hand, a different spectrum associated to $\phi$ has entered recent developments in symplectic topology, starting with the work of Hutchings on the Calabi invariant. The \emph{mean action spectrum} is defined as \[\specm(\phi) = \bigcup_{k \geq 1} \frac{1}{k} \spec(\phi^k) \subset \R.\] It is natural to ask what kind of properties this subset possesses. Though the action spectrum and the mean action spectrum are closely related, an intriguing question of Hutchings suggests that the topology of the mean action spectrum has a strikingly different structure:

\begin{question}[\cite{hutchings2025reeb}]\label{q:H}
Is the closure of the mean action spectrum for an area-preserving map of a surface with zero flux a segment?
\end{question}

This question is positively answered for all autonomous Hamiltonian flows on the standard disk $\mathbb{D} = \lbrace (x,y) \hspace{1 mm} | \hspace{1 mm} x^2 + y^2 \le 1 \rbrace $ in \cite{bechara2023asymptotic} and is also known to hold for monotone twist maps of the annulus, as we explain below.

In the present work, we make progress on Question \ref{q:H} for all symplectic maps of the disk. We first recall the following notions.
\begin{definition}
An area-preserving map $\phi:\mathbb{D} \to \mathbb{D}$ is called a \emph{pseudo-rotation} if the cardinality of the set of periodic orbits of $\phi$ is one. The map $\phi$ is called \emph{periodic} if there exists a natural number $n$ such that $\phi^n = \mathrm{Id}$.
\end{definition}

We partially answer Question \ref{q:H} as follows. 

\begin{theorem}
    If $\phi: \mathbb{D} \to \mathbb{D}$ is a symplectic map and $\widetilde{\varphi}$ is an isotopy from $id$ to $\phi$, then either $\overline{\specm(\varphitilde)}$ contains a positive length segment, or $\overline{\specm(\varphitilde)}$ is a single point. If $\overline{\specm(\varphitilde)}$ is a single point, the map is either a pseudo-rotation or a periodic map.
\end{theorem}

In fact we prove the following quantitatively refined version of this statement.


\begin{theorem}\label{mainthm:meanacion}
Let $\phi:\mathbb{D} \to \mathbb{D}$ be an area-preserving diffeomorphism with isotopy $\phi_{t}$.  
If $\phi$ is neither a pseudo-rotation nor a periodic map, then there exists a segment $[a,b]$ of positive length and positive integers $k$ and $s_{0}$ such that for any subinterval $[c,d] \subset [a,b]$ with $d-c \geq \frac{3}{k2^{s}}$ and for all $s\geq s_{0}$, there exists a periodic point $x_{s}$ of period $k2^{s}$ with 
$$
a_{\phi_{t}}^{\infty}(x_{s}) \in [c,d].
$$  
\end{theorem}

In \cite{bechara2023asymptotic}, the authors showed that if a map is a pseudo-rotation that is conjugate to a rotation on the boundary, then the integral of the action function with respect to any invariant measure equals its rotation number. Theorem \ref{mainthm:meanacion} provides the converse: if a map is neither periodic nor a pseudo-rotation, then there exist two invariant probability measures whose integrals of the action function differ. Indeed, one can take two periodic points of the map with distinct mean actions and use the associated Dirac probability measures supported on these periodic points.

The answer to Question \ref{q:H} for monotone twist maps of an open annulus is known to be positive; see \cite[Definiton 1.1.1]{siburg-princ-least-act-04}. The asymptotic mean action function $a^{\infty}_{\phi,\lambda}$ for a monotone twist map $\phi$ is basically the minimal action in the sense of \cite[Definition 1.3.6]{siburg-princ-least-act-04}. The latter is shown to be strictly convex and continuous for monotone twist maps in \cite[Theorem 1.3.7]{siburg-princ-least-act-04}. On the other hand, every rational number in the twist interval of $\phi$ is equal to the minimal action of a periodic point as shown in \cite[Theorem 1.3.4]{siburg-princ-least-act-04}. Therefore, the closure of the mean action spectrum of monotone twist maps is a connected interval. The billiard maps of strictly convex domains with smooth boundaries are monotone twist maps; see \cite[\S1.2]{siburg-princ-least-act-04}. The mean action of a periodic point of a billiard map corresponds to the length of periodic trajectory of the billiard table divided by the number of times it hits the boundary of the table.

Consequently, one can prove a result similar to Theorem \ref{mainthm:meanacion} for any smooth symplectic diffeomorphism $\phi$ of any surface which admits a \emph{general elliptic} fixed point; see \cite[Definition 4.1.9]{siburg-princ-least-act-04}. Indeed, it follows from Birkhoff's normal form \cite[Theorem 4.1.8]{siburg-princ-least-act-04} that $\phi$ can be conjugated to another symplectic map $\psi$ which looks like a monotone twist map near the elliptic fixed point. So, the result follows from the above discussion; see for instance the proof of \cite[Theorem 4.1.11]{siburg-princ-least-act-04}.



\section*{AI disclosure}
AI was not used in the production of this paper.

\section*{Acknowledgments}
We thank Michael Hutchings for inspiring lectures and encouraging conversations.

D.C-G thanks the NSF for their support under agreements DMS- 2227372 and DMS-2238091.  

H.A. was supported by National Key R\&D Program of China Grant No. 2023YFA1010500 and NSFC Grants No. 12301081. 

The work of A.P. is partially supported by the DFG Walter Benjamin Fellowship, Projektnummer 518128580, and  by the National Science Foundation
under Grant No. DMS-2424441.

E.S. was partially supported by an NSERC Discovery grant, by the Courtois chair in fundamental research, by a Sloan research fellowship, by an FRQNT teams grant, and, at the IAS, by the National Science Foundation under Grant No. DMS-1926686 during the academic year 2024-2025 as Von Neumann fellow and by an anonymous donor during the academic year 2026-2027. 

H.A. and E.S. thank the Institute for Advanced Study for its wonderful research atmosphere during the Spring semester of 2025 when their work on this project began.

\section{Preliminaries}
\subsection{Circle diffeomorphisms and their extensions}
Let $\{\phi_t\}_{t \in [0,1]}$ be an isotopy in the set of orientation preserving homeomorphisms $\mathrm{Homeo}^{+}(S^{1})$ 
from the identity to $\phi_{1} = \phi$. Choose a lift 
$\widetilde{\phi}_t : \mathbb{R} \to \mathbb{R}$ of the isotopy $\phi_{t}$ starting with $\widetilde{\phi}_0 = \mathrm{id}$ ending with $\widetilde{\phi}_{1} =\widetilde{\phi}$. The \emph{rotation number} of $\widetilde{\phi}$ is then defined by
$$
\rho(\tilde{\phi})
= \lim_{n \to \infty} \frac{\widetilde{\phi}^{n}(x) - x}{n} \in \mathbb{R};
$$
this limit exists and is independent of $x$. If the lifted isotopy is changed by an integer shift, $\widetilde{\phi}'_{t} = T^{k}\circ \widetilde{\phi}_{t}$ for some $k \in \mathbb{Z}$ where $T(x) = x + 1$, then $\rho(\widetilde{\phi}_{1}') = \rho(\widetilde{\phi}_{1}) + k$. Hence, the class $[\rho(\widetilde{\phi}_{1})] \in \R/\Z$ is well-defined and we denote it by $\rho(\phi).$ This is the classical rotation number. For $\alpha \in \mathbb{R}$, consider the rotation $R_{\alpha} : S^{1} \to S^{1}$ defined by
$R_{\alpha}(x) = x + \alpha \in \R/\Z$ and its isotopy $R_{\alpha,t}(x) = x + t\alpha$ for $t \in [0,1]$. A lift of this isotopy will be $\widetilde{R}_{\alpha,t}(x) = x + t\alpha,$ which also satisfies $\widetilde{R}_{\alpha,t}(x+1) = \widetilde{R}_{\alpha,t}(x) + 1$. The family $\{\rho(\widetilde{\phi}_1 \circ \widetilde{R}_{\alpha})\}_{\alpha \in \mathbb{R}}$ 
is a path in $\mathrm{Homeo}^{+}(\mathbb{R})$ lifting 
$\{\phi\circ R_{\alpha}\}$ in $\mathrm{Diff}^{+}(S^{1})$. The following classical lemma is about the continuity and the monotonicity of the rotation number.
\begin{lemma}\label{lem:rotcont}
The function $\alpha \mapsto \rho(\widetilde{\phi}_1 \circ \widetilde{R}_{\alpha})$ is continuous and non-decreasing. Moreover, if $\rho(\phi) \in \mathbb{R}\backslash \mathbb{Q}$ then the map is increasing.
\end{lemma}
A similar notion of rotation number can be defined for area-preserving maps of the annulus $A = S^{1}\times [0,1]$. Let $\phi:A\to A$ be an area-preserving homeomorphism. Let $\phi_{t}$ be an area-preserving isotopy with $\phi_{0}=id$ and $\phi_{1}=\phi$. We lift $\phi_{t}$ to the universal cover $\widetilde{A} = \mathbb{R}\times [0,1]$ and denote it by $\widetilde{\phi}_{t}$. The rotation number of $\widetilde{\phi}_{t}$ at a point $z\in A$ is defined by
$$\rho_{\widetilde{\phi}_{t}}(z):=\lim_{n\to \infty}\frac{\pi_{1}\big({\widetilde{\phi}_{1}^{n}(z)}\big)-\pi_{1}(z)}{n},$$
if the limit exists. Here, $\pi_{1}:\mathbb{R}\times [0,1]\to \mathbb{R}$ is the projection to the first coordinate. The class $[\rho_{\widetilde{\phi}_{t}}(z)] \in \R/\Z$ is well-defined and does not depend on the lift of the isotopy and it is denoted by $\rho_{\phi_{t}}(z)$.
\begin{remark}\label{rem:franks_const_rot_pseudo_rot}
In the annulus case, the rotation number might not be well-defined at every point.  Franks \cite{franks1988generalizations} proved that if an area-preserving map $\phi:S^{1}\times [0,1]\to S^{1}\times [0,1]$ does not have a periodic orbit, then $\rho_{\phi_{t}}$ is well-defined and constant. Moreover, the converse is also true; if $\rho_{\phi_{t}}$ is well-defined and equal to an irrational constant, then the map is either periodic or a pseudo-rotation.
\end{remark}

\subsection{Maps with zero topological entropy and their properties}\label{sec:positive_top_entrop}
Topological entropy is a numerical invariant that measures the \emph{complexity} of a dynamical system in a topological sense. 
Intuitively, it quantifies the exponential growth rate of distinguishable orbit segments as time evolves; the higher the entropy is the more ``unpredictable" or ``chaotic” the system is. It is invariant under topological conjugacy, making it a fundamental tool for detecting and classifying complex dynamics, such as the presence of horseshoes and symbolic dynamics.
\begin{definition}\label{def:top_entropy}
Let $(X,d)$ be a compact metric space and let $f\colon X\to X$ be a continuous map.
For every $n\in\mathbb{N}$ define the following metric:
$$d_n(x,y)\ :=\ \max_{0\le k<n} d\bigl(f^k(x),f^k(y)\bigr).$$
Let $\epsilon >0$ be any number. A subset $E\subset X$ is called \emph{$(n,\epsilon)$-separated} if $d_n(x,y)\ge \epsilon$ for every pair of distinct points $x,y\in E$. Let,
$$s(n,\epsilon)\ :=\ \max\bigl\{|E|:\ E\subset X \text{ is }(n,\epsilon)\text{-separated}\bigr\}.$$
The \emph{topological entropy} of $f$ is defined by,
$$h_{\mathrm{top}}(f)\ :=\ 
\lim_{\epsilon\to 0}\,
\limsup_{n\to\infty}\,
\frac{1}{n}\log s(n,\epsilon).$$
\end{definition}
\begin{definition}[Two-sided full shift on two symbols]
Let $\Sigma_2 := \{0,1\}^{\mathbb{Z}}$ be the space of bi-infinite sequences of $0$'s and $1$'s:
$$x = (\dots, x_{-1}, x_0, x_1, \dots); x_i \in \{0,1\}.$$
Endow $\Sigma_2$ with the product topology where $\{0,1\}$ is discrete. 
The \emph{two-sided full shift map} $\sigma\colon \Sigma_2 \to \Sigma_2$ is defined by $(\sigma(x))_i := x_{i+1}$ for $i\in \mathbb{Z}.$
\end{definition}
\begin{remark}
The map $\sigma$ is a homeomorphism of $\Sigma_2$, and its dynamics is highly chaotic: every finite sequence of $0$'s and $1$'s appears somewhere in the orbit of a dense set of points. 
Its topological entropy satisfies $h_{\mathrm{top}}(\sigma) = \log 2$.
\end{remark}
The following result is an implication of Katok's work \cite{katok1980lyapunov}.
\begin{theorem}\label{horseshoe}
Let $f\colon \Sigma\to \Sigma$ be a $C^{1+\alpha}$-diffeomorphism for some $\alpha>0$ of a compact surface $\Sigma$. If the topological entropy of $f$ is positive then $f$ has an iterate $f^N$ that admits a compact invariant set $\Lambda$ on which $f^N|_\Lambda$ is topologically conjugate to a full shift on finitely many symbols.
\end{theorem}
\begin{remark}
Denote $f^{N}$ by $\widehat{f}$. Then, the above theorem implies that for all $n$ the number of hyperbolic fixed points of $\widehat{f}^{n}$ are at least $k^{n}$ for some $k >1$.
\end{remark}

\subsection{Properties of the action spectrum of area-preserving maps}
Let $\phi$ be an area-preserving smooth diffeomorphism of the disc $\mathbb{D}$, isotopic to the identity via an area-preserving isotopy $\{\phi_{t}\}_{t\in [0,1]}$. Let $\lambda$ be a primitive of the standard symplectic form \[\omega=\frac{1}{\pi}\,dx\wedge dy = \frac{r}{\pi}\,dr\wedge d\theta.\] The \textit{action function} of $\phi_{t}$ with respect to $\lambda$ is the unique smooth function $a_{\phi_{t},\lambda}:\mathbb{D}\to \mathbb{R}$ satisfying,
\begin{enumerate}
\item $da_{\phi_{t},\lambda}=\phi^{*}\lambda-\lambda,$
\item $a_{\phi_{t},\lambda}(x)=\int_{\{\phi_{t}(x)\}}\lambda$ for every $x\in \partial\mathbb{D}.$
\end{enumerate}
We have the following identity \cite{abbondandolo2018sharp, Bechara-23-ann-du-quebec} where we use the convention $\omega(X_{t}, -) = -dH_{t}$ for Hamilton's equation; see also Remark \ref{rem:compare_two_action_defs}.
\begin{proposition}\label{prop:rel_two_actions}
Let $\phi_{t}$ be an area-preserving isotopy with time-one map $\phi$ and $H_{t}$ be a generating Hamiltonian of $\phi_{t}$ normalized so that it vanishes on the boundary $\partial \D$. Then, for every $x\in \D$,
\begin{equation}\label{eq:action}
a_{\phi_{t},\lambda}(x) = \int_{\{\phi_{t}(x)\}}\lambda
- \int_{0}^{1} H_t\bigl(\phi_t(x)\bigr)dt.
\end{equation}
\end{proposition}
The following describes how the action changes under varying the isotopy or the primitive.
\begin{lemma}\label{lem:actionfor}
Let $\phi_{t}$ and $\psi_{t}$ be area-preserving isotopies from $id$ to some $\phi$ and $\psi$ respectively. Let $\lambda$ be a primitive of $\omega$ and let $u$ be a smooth real function on $\mathbb{D}$. Then, 
\begin{enumerate}\label{lem:action}
    \item $a_{\phi_{t},\,\lambda + du} 
    = a_{\phi_{t},\lambda} + u \circ \phi - u.$
    \item $a_{\psi_{t}\circ \phi_{t},\lambda} 
  = a_{\psi_{t},\lambda} \circ \phi_{1} + a_{\phi_{t},\lambda} 
    = a_{\psi_{t},\lambda} + a_{\phi_{t},\psi_{1}^{*}\lambda}.$ 
    \item $a_{\phi_{t}^{-1},\lambda} 
    = -\,a_{\phi_{t},\lambda} \circ \phi^{-1} 
    = -\,a_{\phi_{t},(\phi_{1}^{-1})^*\lambda}.$
\end{enumerate}
\end{lemma}
The second property implies that the quantity $a^{\infty}_{\phi_{t}}(x):=a_{\phi_{t}^{n},\lambda}(x)/n$ is independent of the choice of the primitive $\lambda$ for any periodic point $x$ of period $n$ for $\phi$. This allows us to define the \textit{mean action spectrum},
$$\mathcal{S}_{\mathrm{mean}}(\phi_{t})=\left\{a^{\infty}_{\phi_{t}}(x) \,\middle|\, x \text{ is a periodic point of } \phi \right\}.$$
Since $a^{\infty}_{\phi_{t}}$ is independent of the choice of the primitive, for convenience in computations we choose the primitive $\lambda= r^{2}d\theta/2\pi$ of the area form to define the action function. This choice is useful because when the map is radial on some region the action has the same value along the entire circle in that region. From Lemma \ref{lem:action}, we have the following. 
\begin{lemma}\label{lem:congact}
Let $\phi_{t}$ be an area-preserving isotopy from $id$ to $\phi$ and $\psi$ be another area-preserving diffeomorphism. Then, for every periodic point $x$ of $\phi,$ we have
\begin{enumerate}
    \item $a^{\infty}_{\phi_{t}}(x)=a^{\infty}_{\psi\circ \phi_{t}\circ \psi^{-1}}(x),$
    \item $a^{\infty}_{\phi_{t}^{n}}(x)=n a^{\infty}_{\phi_{t}}(x),$
    \item $a^{\infty}_{\phi_{t}\circ R_{2\pi nt}}(x)=a^{\infty}_{\phi_{t}}(x)+n,$  
\end{enumerate}
where $R_{2\pi nt}(r,\theta):=(r,\theta+2\pi nt)$ is the standard rotation.
 \end{lemma}
The proof of Lemma \ref{lem:congact} is an immediate application of Lemma \ref{lem:action}. 
\subsection{Action spectrum and winding number.}
We first define the winding number for an isotopy $\phi_{t}$ of area-preserving diffeomorphisms from $id$ to $\phi.$
\begin{definition}
The \textit{winding number} $W_{\phi_{t}}(x, y)$ of a pair of distinct  points $x,y\in \D$ is,
$$W_{\phi_{t}}(x, y):=(\theta(1)-\theta(0))/2\pi,$$
where $\theta: [0,1] \rightarrow \R$ is a continuous function such that
$\vert\phi_t(y)-\phi_t(x)\vert e^{i\theta(t)}= \phi_t(y)-\phi_t(x)$ for all $t \in [0,1].$
\end{definition}
We define the asymptotic winding number by \[W^{\infty}_{\phi_{t}}(x, y):=\lim_{n \to \infty} W_{\phi_{t}^{n}}(x, y)/n\] provided that the limit exists. The following result proven in the Ph.D. thesis of Pu and later generalized by David-Bechara for every point and not just for periodic points; see \cite[Theorem 4.1]{bechara2023asymptotic} and also \cite[Theorem 4.3.2]{senior2023asymptotic}.
\begin{theorem}\label{thm:actwind}
For almost every point $x \in \D$, including periodic points of $\phi$, one has
$$a^{\infty}_{\phi_{t}}(x)=\int_{\mathbb{D}}W^{\infty}_{\phi_{t}}(y,x)\omega(y).$$
\end{theorem}
\begin{remark}
The winding number $W_{\phi_{t}}(x, y)$ is uniformly bounded and well-defined for every $x \neq y$. This was discussed in detail in \cite{senior2023asymptotic}.
\end{remark}

\subsection{Compactly supported extensions}\label{section:approx}
We assume that $H\in C^\infty\big(\R/\Z\times\D\big)$ satisfies
\begin{equation}\label{eq1}
H_t(z) = 0;  \forall t\in\R/\Z, z\in\partial\D.
\end{equation}
To extend $H$ to a compactly supported Hamiltonian, on a slightly larger disk, we take two different approaches depending on the rationality of the rotation number $\rho(\phi_t)$ on the boundary $\partial \D$ of the isotopy $\phi_{t}$ generated by $H$.
\subsubsection{Rational rotation number}\label{sec:comp_supp_ext_rational}
As we care about the mean actions of periodic points, up to passing to an iteration of $H$ and switching to another appropriate Hamiltonian generating the same time-one map we can assume that the rotation number on the boundary is zero. To construct a compactly supported extension of $H$, we use a sequence of smooth functions $\rho_{n}:\R\to\R$ satisfying
\begin{enumerate}
    \item $\rho_n(r)=r$ on $r\leq0$,
    \item $\rho_n(r)=0$ on $r\geq1/n$,
    \item $\rho_n\geq0$ on $r\geq0$,
    \item $\vert\rho_{n}'\vert_{C^{0}}\leq 1$ and $\vert\rho_n\vert_{C^0}\leq 1/n$ on $r \geq 0$. 
\end{enumerate} 
We fix an arbitrary smooth extension of $H$ to the entire plane and continue to denote it by $H.$ Note that this extension exists as $H$ is a smooth function on $\R/\Z \times \D.$ In polar coordinates this gives a function $H(t,r,\theta)$ on $\mathbb{R}/\mathbb{Z} \times [0,\infty) \times \mathbb{R}/\Z$. Define
\begin{equation}\label{D:H_n}
H_n(t,1+r,\theta) := H\bigl(t,1+\rho_n(r),\theta\bigr); n \in \mathbb{N}.
\end{equation}
In \cite{bramham2024spectral}, the authors proved that the sequence $H_n$ satisfies the following properties:
\begin{enumerate}
\item (Extension) $H_n$ restricted to $[0,1]\times\D$ is equal to $H$. \item (Convergence) $\lim_{n\to\infty}\vert H_n\vert_{C^0(\R^2\backslash\D)}=0$.  
\end{enumerate}
The following property proved in \cite{bramham2024spectral} holds for the periodic orbits of the flow $\phi_{n,t}$ generated by $H_n$ in the extended region.
\begin{lemma}\label{L:wind}
Suppose that the generating Hamiltonian $H$ has an autonomous Hamiltonian vector field on the boundary $\partial\D$. Let $H_n$ be given by \eqref{D:H_n} and suppose that the rotation number on the boundary satisfies  $|\rho(\phi_{n,t})|<1.$ Then, for large enough $n$ the winding number of every one-periodic orbit $x:[0,1]\to\R^2 \backslash \D$ of $\phi_{n,t}$ around the origin satisfies $|\wind(x)|<1$. 
\end{lemma} 

\begin{remark}\label{ext:smallact}
As explained in the proof of \cite[Theorem 4]{bramham2024spectral}, Lemma \ref{L:wind} implies that for arbitrary large $n$, the actions of fixed points of $\phi_{n,1}$ in $\mathbb{R}^{2}\backslash \mathbb{D}$ will be arbitrarily close to $0$. So, if $\rho(\phi_{n,t}) = 0$, so will be $\rho(\phi_{n,t}^{k}) = 0$. Hence, for a fixed number $k\geq 1$, applying the above discussion to the $k$-th iteration of $H_{n}$, one proves that that the mean actions of periodic points of periods at most $k$ of $\phi_{n,1}$ get arbitrary close to zero as $n \to \infty$. See the proof of \cite[Theorem 4]{bramham2024spectral} for more details.
\end{remark}

\begin{remark}
Later we introduce another extension method for the case of irrational rotation number on the boundary that differs from the rational case above. The reason is the following. For any area preserving diffeomorphism $\phi$ one can always find a Hamiltonian $H$ generating a flow $\phi_{t}$ with $\phi_{1} = \phi$ and rotation number $\vert\rho(\phi_{t})\vert < 1$. If the rotation number is irrational, one cannot find such a flow $\phi_{t}$ with small rotation number so that all of its iterates $\phi_{t}^{k}$ also have small rotation numbers on the boundary. Therefore, the mean actions of periodic orbits of high periods are not necessarily close to zero; this was necessary for instance in the proof of Theorem \ref{thm:ratcase}.
\end{remark}
\subsubsection{Irrational rotation number}\label{sec:comp_supp_ext_irrational}
In this situation, we perform a $C^{\infty}$-small perturbation that creates no periodic points of periods less than a given large integer $k$. After the perturbation, the map is area-preserving
conjugate to a map that is a rotation on the boundary, and we shall work with that model. We start with the following definition of irrational numbers.
\begin{definition}\label{def:dioph_num}
An irrational number $x$ is called \emph{Diophantine} if there exist 
constants $\sigma \ge 0$ and $\gamma > 0$ such that,
$$\left| x - \frac{p}{q} \right| \ge \frac{\gamma}{q^{2+\sigma}}$$
for all integers $p$ and $q \neq 0$. 
The set of Diophantine numbers has full measure in $\mathbb{R}$.
\end{definition}
In the following construction, we will use a deep theorem of Herman and Yoccoz \cite{herman1979conjugaison, yoccoz1984conjugaison} 
on circle diffeomorphisms with Diophantine rotation number being conjugates of rigid rotations.
\begin{theorem} \label{thm:rotonboundary} Let $\phi:\mathbb{D}\to\mathbb{D}$ be an area-preserving diffeomorphism such that the rotation number of $\phi$ restricted to the boundary is irrational. Then, for every $n\geq 1$ there exists a smooth approximation $\phi_{n}$ of $\phi$ such that, 
\begin{enumerate}
    \item $|a_{\phi}-a_{\phi_{n}}|_{C^{0}}\leq 1/n$.
    \item For all $k\leq n$, the set of periodic points of period $k$ for $\phi_{n}$ and $\phi$ coincide.
    \item There exists smooth an area-preserving diffeomorphism $\nu_{n}$ such that in polar coordinates $\nu\circ \phi_{n}\circ \nu^{-1}(1,\theta)=(1,\theta+2\pi\theta_{n})$ for an irrational number $\theta_{n}$. 
\end{enumerate}
\end{theorem}
\begin{proof}
Since $\phi$ has irrational rotation number on the boundary there are small positive numbers $\eta, \delta> 0$ so that $\phi^{j}$ shifts every point $x\in \D \backslash \D_{1-\eta}$ by $2\delta$, i.e., one has $d(\phi^{j}(x),x)\geq 2\delta$ for every $j \in \{1,2,\dots, n\}$. Define $\psi_{t}: \D \rightarrow \D$ to be a Hamiltonian isotopy generated by a $C^{1}$-small radial Hamiltonian $h(r)$ that vanishes on $\D_{1-\eta}$ and is linear in $r$ with small slope $-\epsilon$ outside of $\D_{1-\frac{\eta}{2}}$ for some small number $\epsilon > 0$ and $h''(r)\leq 0$ everywhere Then, it is clear that $(\phi \psi)^{j}$ still $\delta$-shifts all the points $x\in \D\backslash \D_{1-\eta}$ for every $j\in \{1,\dots, n\}$ if $\epsilon$ is small enough. Hence, one creates no new $j$-periodic point for $\phi\psi$ compared to $\phi$ for all $1\leq j\leq n$. On the other hand for $\lambda = r^{2}d\theta/2\pi$,
$$da_{\psi_{t},\lambda}=\psi_{1}^{*}(r^{2}\d\theta/2\pi) - r^{2}\d\theta/2\pi =r^{2}h'(r)dr/2\pi.$$
Given that $\vert h\vert_{C^{1}} \leq \epsilon$, this implies that $\lvert a_{\psi_{t},\lambda}\rvert\leq \epsilon.$ Together with property $(2)$ of Lemma \ref{lem:actionfor} this yields
$\lvert a_{\phi_{t} \psi_{t}}-a_{\phi_{t}}\rvert\leq \epsilon.$ Hence, the map $\phi \psi$, and in fact all the maps $\phi \psi_{t}$ for $t\in [0,1]$, satisfy the items $(1)$ and $(2)$ of the theorem. The restriction of $\phi \psi_{t}$ to the boundary of the disc takes the form $\phi \psi_{t}|_{\partial \mathbb{D}} = \phi|_{\partial \mathbb{D}} \circ R_{-\epsilon t}$ where $R_{a}$ is the rotation by $2\pi a$. By Lemma \ref{lem:rotcont} we can choose a sufficiently small $t$ such that the rotation number of $\phi \psi_{t}$ on the boundary $\partial \mathbb{D}$ is Diophantine; see Definition \ref{def:dioph_num}. By a result of Herman and Yoccoz \cite{herman1979conjugaison, yoccoz1984conjugaison}, for such $t,$ there exists a smooth circle diffeomorphism $\nu$ such that 
$\nu \circ (\phi \psi_{t})|_{\partial \mathbb{D}} \circ \nu^{-1}$ is a rigid rotation. To complete the proof of the theorem, we need to extend $\nu$ to an area-preserving diffeomorphism of $\mathbb{D}$. There are several ways to extend a circle diffeomorphism to a smooth area-preserving diffeomorphism of the disk, e.g., one can use the construction of Banyaga \cite{banyaga1974formes}.
\end{proof}
We now state a result that explains how to extend maps that are rotations on the boundary of 
$\mathbb{D}$ to a larger disk, in such a way that the extended map is a rotation near the boundary and has no periodic orbits of period less than or equal to $n$ for a given $n$. In the following $\D_{1 + r}$ is the disk of radius $1 + r$. 
\begin{theorem}\label{thm:extensionrot}
Assume that the map $\phi:\mathbb{D}\to \mathbb{D}$ is an area-preserving diffeomorphism and irrational rigid rotation on the boundary that is  $\phi(1,\theta)=(1,\theta+2\pi\theta_{0})$ for some $\theta_{0}\in \mathbb{R}\backslash \mathbb{Q}$. Then, for all $n\in \mathbb{N}$, there exists a smooth extension $\phi_{n}$ of $\phi$ to $\D_{1+1/n}$ such that,
    \begin{enumerate}
    \item $\phi_{n}(r,\theta)=(r,\theta+2\pi\theta_{0})$ near the boundary $\partial \D_{1 + 1/n}$.
    \item The map  $\phi_{n}$ has no new periodic points of periods less than $n$ compared to $\phi$.
    \item There exist isotopies $\phi_{t}, \phi_{n,t}$ with time one maps $\phi, \phi_{n}$ on $\D, \D_{1 + 1/n}$ so that \[a^{\infty}_{\phi_{n,t}}(x) = a^{\infty}_{\phi_{t}}(x) + a_{n}\] for all $x\in \D$ where $a_{n}$ is a constant satisfying $|a_{n}|\leq 3\theta_{0}/n$.
    \item Finally one has $|a_{\phi_{n,t}}(x)-\mathrm{area}(\D_{1+1/n})\theta_{0}|\leq 1/n$ for every $x\in \D_{1+1/n}\backslash \mathbb{D}.$
    \end{enumerate}
\end{theorem}
The construction of $\phi_n$ is very similar to the proof of Proposition 4.1 in \cite{pirnapasov2024hutchings}, so we omit the full details. The idea is to first extend it by the rotation $2\pi \theta_{0}$. The extended map will be a $C^\infty$ area-preserving diffeomorphism except on the boundary of $\mathbb{D}$. Then, one can use a suitable generating function and the Whitney approximation theorem \cite{whitney1934differentiable} to obtain the desired smooth approximation. In addition, this construction helps to control the action function of $\phi_n$ in the extended region. The following were proven by Alves and Meiwes \cite{alves2021braid}.
\begin{lemma}\label{lem:specialham}
If $\phi:\mathbb{D}\to\mathbb{D}$ is an area-preserving diffeomorphism such that $\phi(r,\theta)=(r,\theta+2\pi\theta_{0})$ for $r\geq 1-2\epsilon$ and some $\epsilon > 0$, then, there is a Hamiltonian function $H:S^{1}\times\mathbb{D}\to\mathbb{R}$ generating $\phi$ as its time one map such that $H_{t}(r,\theta)=\theta_{0}(r^2-1)$ for all $r\geq 1-\epsilon.$
\end{lemma}
\begin{remark}
Near the boundary, the map can be expressed as composition of a rotation and a compactly supported map. Applying the construction in Lemma 3.5 of \cite{alves2021braid}, we obtain an appropriate isotopy extending this rotation.
\end{remark}
Now we extend the Hamiltonian functions from Lemma \ref{lem:specialham} as compactly supported Hamiltonians after shifting it by a small constant, while controlling the action of fixed points in the extended region. To do this, we use auxiliary cut-off functions $f_{\beta, \theta_{0}}: [1, 1+\beta] \to \mathbb{R}$ for $\beta, \theta_{0}> 0$ satisfying:
\begin{enumerate}
    \item $f_{\beta,\theta_{0}}(r) = \theta_0$ for $r$ sufficiently close to $1$,
    \item $f_{\beta,\theta_{0}}(r) = 0$ for $r$ sufficiently close to $1 + \beta$,
    \item $f_{\beta,\theta_{0}}(r)$ is strictly decreasing when $f\neq\theta_0$ or $f\neq 0$, 
    \item $f_{\beta,\theta_{0}}$ varies smoothly as a function of $\beta$.
\end{enumerate}
Assume $\phi:\mathbb{D}\to\mathbb{D}$ is an area-preserving map with $\phi(r,\theta)=(r,\theta+2\pi\theta_{0})$ on the region $r\geq 1-2\epsilon$ for  some $\epsilon > 0$. Let $H$ generate the Hamiltonian isotopy $\phi_{t}$ with time one map $\phi$ given by Lemma \ref{lem:specialham}, in particular we have $\phi_{t}(r, \theta) = (r, \theta + 2\pi \theta_{0}t)$ near the boundary $\partial \D$. Consider the extended isotopy,
$$\phi_{\beta,t}(r,\theta):=
\begin{cases}
\phi_{t}(r,\theta), & r \leq 1,\\[6pt]
(r, \theta + 2\pi f_{\beta, \theta_{0}}(r)t), & r \geq  1.
\end{cases}$$
to the slightly larger disk $\D_{1 + \beta}$ and set $\phi_{\beta} = \phi_{\beta, 1}$ which is a compactly supported Hamiltonian diffeomorphism of $\D_{1 + \beta}$. Recall that the action function $a_{\phi_{\beta},\lambda}$ is well-defined and does not depend on the choice of the isotopy generating $\phi_{\beta}$; here $\lambda = r^{2}d\theta/2\pi$. In the following, for brevity, we set $f = f_{\beta, \theta_{0}}$.

\begin{lemma}\label{lem:extcompact}  
The action function $a_{\phi_{\beta}, \lambda}: \D_{1 + \beta} \rightarrow \R$ satisfies the following properties:
\begin{enumerate}
    \item For all $1 \leq r \leq 1 + \beta$,
    $$a_{\phi_{\beta}, \lambda}(r,\theta)=r^{2}f(r)+\int^{1+\beta}_{r}2sf(s)ds.$$ 
    \item $a_{\phi_{\beta}, \lambda}(x) = a_{\phi_{t}, \lambda}(x) + a_{\beta}$ for every $x\in \mathbb{D}$ where \[a_{\beta} = \int^{1+\beta}_{1}2sf(s)ds > 0, \]
    $$|a_{\beta}|=|\int^{1+\beta}_{1}2sf(s)ds|\leq \theta_{0} \beta(2+\beta).$$
    \item $a_{\phi_{\beta},\lambda}(r,\theta)$ is decreasing as a function of $r$ on the region $1 \leq r \leq 1+\beta$.
    \item A point $(r,\theta)\in \D_{1 + \beta}$ is a periodic point of $\phi_{\beta}$ if $kf(r)$ is a rational number. Moreover, if $kf(r_{1})=n+1$ and $kf(r_{2})=n$ and $r_{1} < r_{2}$, then for all $\theta_{1}, \theta_{2}$ we have, 
    $$1/k\leq a_{\phi_{\beta}, \lambda}(r_{1},\theta_{1}) - a_{\phi_{\beta}, \lambda}(r_{2},\theta_{2})\leq (1+\beta)^{2}/k$$ 
\end{enumerate}
\end{lemma}
\begin{proof}
By definition the action function $a_{\phi_{\beta}, \lambda}$ satisfies the following:
\begin{enumerate}
\item $da_{\phi_{\beta}, \lambda}=\phi_{\beta}^{*}(r^{2}d\theta/2\pi) -r^{2}d\theta/2\pi$,
\item $a_{\phi_{\beta}, \lambda}=0$ on the boundary $\partial \mathbb{D}_{1+\beta}$.
\end{enumerate}
This implies that on the region $1 \leq r \leq 1 + \beta$, 
\begin{equation}\label{form:actbeta}
a_{\phi_{\beta}, \lambda}(r,\theta)=-\int_{r}^{1+\beta}s^{2}f'(s)ds =r^{2}f(r)+\int^{1+\beta}_{r}2sf(s)ds.
\end{equation}
Since $\phi_{\beta}$ and $\phi$ coincide on $\mathbb{D}$ the action functions $a_{\phi_{\beta}, \lambda}$ and $a_{\phi_{t}, \lambda}$  differ by a constant. We know that $a_{\phi_{t}, \lambda}(1,\theta)=\theta_{0}$ and $a_{\phi_{\beta}, \lambda}(1,\theta)=\theta_{0}+\int_{1}^{1+\beta}2sf(s)ds$ which imply $a_{\phi_{\beta}, \lambda} = a_{\phi_{t}, \lambda} + a_{\beta}$ where 
$$|a_{\beta}| = |\int_{1}^{1 + \beta}2sf(s)ds|\leq \beta(2+\beta)\theta_{0}.$$
By taking differential of $a_{\phi_{\beta}, \lambda}$ in the radial direction one observes that for $r\in [1,1+\beta]$ one has $\partial_{r}a_{\phi_{\beta}, \lambda}=r^{2}f'(r)\leq 0$. Finally, using the formula \ref{form:actbeta} for all $\theta_{1}, \theta_{2}$ we get,
\begin{equation}\label{eq:diff_actions_rad}
a_{\phi_{\beta}, \lambda}(r_{1},\theta_{1})-a_{\phi_{\beta}, \lambda}(r_{2},\theta_{2})=r_{1}^{2}f(r_{1}) - r_{2}^{2}f(r_{2}) + \int_{r_{1}}^{r_{2}}2sf(s)ds.
\end{equation}
Observe that $n\leq kf(s)\leq n+1$ for every $s\in [r_{1},r_{2}]$ as $f(s)$ is decreasing on the interval $[r_{1}, r_{2}]$ and in particular $f(r_{2}) \leq f(s) \leq f(r_{1})$. This implies,
\begin{equation}
\begin{aligned}
1/k \leq (f(r_{1}) - f(r_{2}))r_{1}^{2} \leq \eqref{eq:diff_actions_rad} \leq (f(r_{1}) - f(r_{2}))r_{2}^{2} \leq (1 + \beta)^{2}/k.
\end{aligned}
\end{equation}
This completes the proof of the Lemma \ref{lem:extcompact}.
 \end{proof}
\begin{remark}
By construction of $\phi_{\beta}$ for the extended region, one can easily show that these maps are Morse-Bott non-degenerate in $\mathbb{D}_{1+\beta}\backslash \mathbb{D}.$
\end{remark}
\begin{remark}\label{rmk:easier_extension}
Here is another possibly easier way to see the extended Hamiltonian $\phi_{\beta,t}$ and prove the properties of its action function. Consider $\phi_{t}$ generated by a Hamiltonian $H$ that is of the form $\theta_{0}(r^{2} - 1)$ near the boundary $\partial \D$. Extend $H$ to a Hamiltonian $H_{\beta} + a_{\beta}$ on $\D_{1 + \beta}$ that stays radial and its slope decreases from $\theta_{0}$ to zero and its value on $\partial \D_{1+\beta}$ is a small constant value $a_{\beta}$ which converges to zero as $\beta \to 0$; see Figure \ref{fig:extended_hamilt}. Let $\phi_{\beta,t}$ be generated by $H_{\beta}$; note that the latter Hamiltonian vanishes near the boundary $\partial \D_{1 + \beta}$ and so we can use it to define the spectral invariant of $\phi_{\beta,t}$. It is now clear that for fixed points $x\in \D$ one has $\mathcal{A}_{H_{\beta}}(x) = \mathcal{A}_{H}(x) - a_{\beta}$. The action of a fixed point $x = (r,\theta)\in \D_{1+\beta}\backslash \D$ of $\phi_{\beta,1}$ equals $-\theta_{0}$ on $r=1$ and increases to zero as $r$ increases. It is also easy to confirm that $\mathcal{A}_{H_{\beta}}(r_{2},\theta_{2}) - \mathcal{A}_{H_{\beta}}(r_{1},\theta_{1})$ for $r_{2} > r_{1}$ is a small positive number provided $r_{1},r_{2}$ are close. These properties translate to the properties stated in Lemma \ref{lem:extcompact} observing that $\mathcal{A}_{H} = -a_{\phi_{t}}$ and $\mathcal{A}_{H_{\beta}} = -a_{\phi_{\beta,t}}$; see Remark \ref{rem:compare_two_action_defs}. It is also easy to show, using continuity and spectrality of spectral invariants, that $c_{+}(\phi_{\beta,t}) = \sigma_{+}(\phi_{t}) - a_{\beta}$ provided $c_{+}(\phi_{t}) > 0$ and $\beta$ is small enough.
\begin{figure}[ht]
  \centering
  \begin{tikzpicture}
\draw (-1,-1.5) -- (-1,1);
\draw[->] (-1,0) -- (3, 0);\node at (3.3, 0) {$r^{2}$};
\draw[line width=0.75pt] (-0.5,-1) -- (0,0);
\draw[dotted][line width=1pt] (-0.75,-1.5) -- (-0.5,-1);
\draw[line width=0.75pt] (1,0.5) -- (2.5,0.5);
\draw[dashed] (-1,0.5) -- (1,0.5); 
\draw[dashed] (1,0) -- (1,0.5); 
\node at (-1.3, 0.5) {\small{$a_{\beta}$}};
\node at (0.1, -0.2) {\tiny{$1$}};\node at (1, -0.2) {\tiny{$(1+\beta)^{2}$}};
\draw[line width=0.75pt] (0,0)to[out=60,in=180,looseness=1.2](1,0.5);
\end{tikzpicture}
\caption{The shifted extended Hamiltonian $H_{\beta} + a_{\beta}$.}
\label{fig:extended_hamilt}
\end{figure}
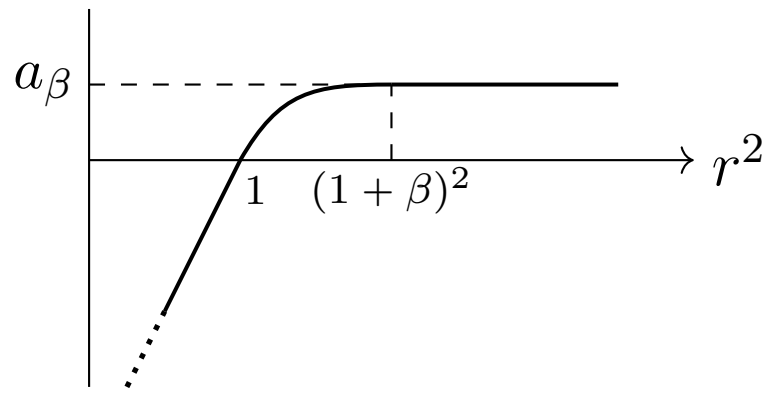
\end{remark}
\subsection{Floer theory background and Smith inequality}\label{sec:floe_theo_back_smith_ineq}
In this section we quickly review Floer theory, focusing on the two dimensional disk with area form $\omega = d\lambda_{st}$ with $\lambda_{st} = r^{2}d\theta/2\pi$ as a Liouville domain, and fix our conventions; for an incomplete list of works on Floer theory in non-compact settings see \cite{floer-hofer-salamon-transversality-95, Viterbo-99, ginzburg-05, lanzat-13, Lanzat-16, alizadeh-atallah-cant}. Then, we state the Smith inequality for Floer homology of compactly supported Hamiltonian diffeomorphisms which is the main tool of the paper.

\subsubsection{Floer homology}\label{sec:floer_cohom}
Let $H: S^{1} \times \D \rightarrow \R$ be a Hamiltonian function compactly supported inside the disk $\D$. One defines its Hamiltonian flow $\varphi_{t}$ via
$$\partial_{t}\varphi_{t} = X_{t}\circ \varphi_{t},\ \omega(X_{t}, )= -dH_{t}$$ for $H_t(x) = H(t,x)$ for all $t\in S^1, x\in \D.$ Given a fixed point $p$ of $\varphi_{1}$ the loops $\varphi_{t}(p), t\in [0,1]$ form the critical points of the symplectic action functional,
\begin{equation}\label{eq:ham_act_funct}
\mathcal{A}_{H}: C^{\infty}(S^{1}, \D) \rightarrow \R;\ \gamma\mapsto \int_{S^{1}} H(t, \gamma(t)) dt - \int_{\D} u^{*}\omega,
\end{equation}
where $u: \D \rightarrow \D$ is a smooth map that coincides with $\gamma$ on the boundary $\partial\D$ under the identification $S^1 \to \partial\D,$ $\theta \mapsto e^{2\pi \theta i}.$ 
\begin{remark}\label{rem:compare_two_action_defs}
If $p\in \D$ is a fixed point of $\varphi_{1}$ and $u$ is a disk with boundary $\{\varphi_{t}(p)\}_{t\in [0,1]}$ then for a Hamiltonian $H$ generating $\varphi_{t}$ and vanishing on $\partial \D$, by Proposition \ref{prop:rel_two_actions} we have $\mathcal{A}_{H}(p) = -a_{\varphi_{t},\lambda}(p)$ where by $\mathcal{A}_{H}(p)$ we mean symplectic action $\mathcal{A}_{H}$ of the loop $\{\varphi_{t}(p)\}_{t\in [0,1 ]}$. Note that, the structure of the mean action spectrum with respect to $\mathcal{A}_{H}$ and $a_{\varphi_{t},\lambda}$ are the same, in particular the closure of one contains a non-trivial interval if and only if the other does. So, we can freely work with each one of them in the subsequent sections.
\end{remark}
Consider a small non-degenerate perturbation of $H$ that is one-homogeneous linear in symplectic radial co-ordinate outside of $\D$ with small positive slope $\epsilon > 0$, i.e., $H = \epsilon r^{2}$ near $\partial \D$, and denote it by $H_{\epsilon}$. The underlying $\Z/2\Z$-vector space of the Floer complex is generated by the critical points of the symplectic action function associated with $H_{\epsilon}$. One picks a family of almost complex complex structures $j_{t}, t\in [0,1],$ compatible with $\omega$ that coincides with the standard complex structure $j_{st}$ at infinity and defines the Floer differential on the Floer complex by counting isolated solutions to Floer's equation:
$$d\gamma_{-} = \sum_{\gamma_{+}} \# \mathcal{M}_{0}(\gamma_{-}, \gamma_{+}, H_{\epsilon},j) \gamma_{+},$$
where $\mathcal{M}_{0}(\gamma_{-}, \gamma_{+}, H_{\epsilon},j)$ is the zero dimensional component of the quotient by the free $\R$-action by translation on the moduli space of smooth maps $u: \R \times S^{1} \rightarrow \C$ satisfying
$$\partial_{s}u + j_{t}(u)(\partial_{t}u - X_{t}(u)) = 0;\ u(\pm \infty, t) = \gamma_{\pm}(t).$$
Floer homology is defined to be the homology of the Floer complex $\CF(H_{\epsilon},j)$ with the differential defined above and is denoted by $\HF(H_{\epsilon})$. The Floer differential decreases the action. The action functional $\mathcal{A}_{H_{\epsilon}}$ defines a filtration on the Floer complex, we denote by $\CF^{\leq t}(H_{\epsilon},j)$ the sub-complex generated by the orbits with action smaller than $t$ and by $\HF^{\leq t}(H_{\epsilon})$ its homology. 

For a compactly supported Hamiltonian $H$ and two non-degenerate perturbations $H_{\epsilon_{1}}, H_{\epsilon_{2}}$ that are one-homogeneous linear at infinity with small positive slopes $\epsilon_{1} \leq \epsilon_{2}$, both smaller than the minimal period of the Reeb flow of $\lambda_{st}$, there is a canonical isomorphism called continuation map,
$$\mathfrak{c}: \HF(H_{\epsilon_{2}}) \rightarrow \HF(H_{\epsilon_{1}}).$$
We define the Floer homology $\HF^{\leq t}(H)$ to be the limit of $\HF^{\leq t}(H_{\epsilon})$ where $\epsilon$ and the perturbation datum inside the disk converge to zero. Indeed, all of these Floer homology groups are isomorphic with $\H_{*}(\D)$ via PSS-isomorphisms which commute with the continuation maps. We denote the image of the degree zero generator of $H_{*}(\D)$, the class of the minimum of a Morse function, by $\min\in \HF_{*}(H_{\epsilon})$.

Given a compactly supported Hamiltonian $H$, on the other hand, one can analogously consider non-degenerate perturbations $H_{-\epsilon}$ that are one-homogeneous linear with small negative irrational slope near the boundary $\partial\D$ and Floer homology groups $\HF_{*}(H_{-\epsilon})$ which are isomorphic with $\H_{*}(\D, \partial \D)$ via some PSS-isomorphisms. In this case, we denote by $\mathrm{max}\in \HF_{*}(H_{-\epsilon})$ the corresponding class to the class of the maximum of a Morse function in degree $2$ of $\H_{*}(\D, \partial\D)$. There are also continuation maps $\mathfrak{c}: \HF(H_{-\epsilon_{1}}) \rightarrow \HF(H_{-\epsilon_{2}})$ for every $0 < \epsilon_{1} \leq \epsilon_{2}$ that commute with the PSS-maps. The filtered chain complexes are $\CF^{\leq t}(H_{-\epsilon})$ consists of periodic orbits with action smaller than $t$.

\subsubsection{Spectral invariants} 
We recall the spectral invariants defined in \cite{viterbo1992symplectic, Schwarz-00, oh-spec-inv-05, Frauenfelder-Schlenk-07}. Given a Hamiltonian function $H$ compactly supported in the disk $\D$ we define 
$$c_{-}(H):= \lim_{\epsilon \to 0}\ \inf\{t: \min\in \mathrm{img}(\HF^{\leq t}(H_{\epsilon}) \rightarrow \HF(H_{\epsilon}))\},$$
where by $\epsilon \to 0$ we mean that the perturbation datum consisting of the slope $\epsilon$ and non-degenerate perturbation inside the disk go to zero. We also define
$$c_{+}(H):= \lim_{\epsilon \to 0}\ \inf\{t: \max\in \mathrm{img}(\HF^{\leq t}(H_{-\epsilon}) \rightarrow \HF(H_{-\epsilon}))\}.$$
It is not hard to show the following relation:
$$c_{+}(H) = - c_{-}(\bar{H}),$$
where $\bar{H}_t = -H_t\circ \phi_{t}$ and $\phi_{t}$ is the Hamiltonian flow of $H$. One can show that the numbers $c_{\pm}$ only depend on the time-one map of the Hamiltonian flow of $H$ as long as $H$ is compactly supported. Hence, for a compactly supported Hamiltonian diffeomorphism $\phi: \D \rightarrow \D$ we define $c_{\pm}(\phi) = c_{\pm}(H)$ where $H$ is a compactly supported Hamiltonian generating $\phi$ as its time-one map. Similarly one can define $c_{+}$ for Hamiltonian functions that are of the form $-\epsilon r^{2}$ near $\partial \D$ for some $\epsilon > 0$. For a Hamiltonian diffeomorphism $\varphi$ that is rotation near the boundary with angle $-2\pi\epsilon$ we define $c_{+}(\varphi):= c_{+}(H_{- \epsilon})$ where $H_{-\epsilon} = - \epsilon r^{2}$ near the boundary and generates $\varphi$ as its time-one map. One also defines $c_{-}(\varphi) := c_{-}(H_{\epsilon})$ where $\varphi$ is Hamiltonian diffeomorphism that is rotation by angle $2\pi \epsilon$ near the boundary, $H_{\epsilon} = \epsilon r^{2}$ near the boundary and generates $\varphi$ as its time-one map.
\begin{definition}\label{D:Hofer norm}
The spectral norm of a compactly supported Hamiltonian diffeomorphism $\phi$ of the disk $\D$ is defined by $\gamma(\phi):=c_+(\phi)-c_-(\phi)$ and the Hofer norm of $\phi$ is $\vert \phi\vert_{\Hof}:=\inf_{H}\vert H\vert_{\Hof}$ where the infimum is taken over all compactly supported $H$ with time-one map $\phi$ and 
$$
\vert H\vert_{\Hof} :=\int_0^1\big(\max_{z\in \D}H(t,z)-\min_{z\in \D}H(t,z)\big) dt
$$  
\end{definition}
Spectral invariants satisfy the following properties; see \cite{viterbo1992symplectic} where the spectral invariants were originally defined using generating functions which later were shown \cite{milinkovic-oh-97} to coincide with the above invariants.

\begin{theorem}\label{T:Vit 1}
For any compactly supported Hamiltonian diffeomorphism $\phi$ generated as the time-one map of $H\in C_{c}^\infty(S^{1}\times\D)$ the spectral invariants have the following properties: 
\begin{enumerate}
 \item $c_+(\phi^{-1})=-c_-(\phi)$. 
 \item $c_-(\phi)\leq 0\leq c_+(\phi)$. 
 \item $c_+(\phi\psi)\leq c_+(\phi)+c_+(\psi)$.
 \item $\mathrm{area}(U) \leq 2\gamma(\phi)$ if $\phi(U)\cap U=\emptyset$; where $\gamma$ is the spectral norm.
 \item $\gamma(\phi)\leq\mathrm{area}(\D)$; see \cite[Theorem 1.1]{alizadeh-atallah-cant-ellipsoid-24}.
 \item $\gamma(\phi)\leq \vert \phi\vert_{\Hof}$, where the latter is the Hofer nom.
 \item $\vert c_{\pm}(\phi_{1})-c_{\pm}(\phi_{2})\vert \leq \vert H_{1} - H_{2}\vert_{\Hof}$ where $H_{1}, H_{2}\in C_{c}^\infty(S^{1}\times\D)$ generate $\phi_{1}, \phi_{2}$.
\end{enumerate}
\end{theorem}
In the following lemma the area of the disk is equal to $1$. 
\begin{lemma}\label{lem:pos_spec_inv}
For a compactly supported Hamiltonian diffeomorphism $\phi$ of the disk,
$$-1 \leq c_{-}(\phi) \leq 0 \leq c_{+}(\phi) \leq 1.$$
\end{lemma}
\begin{proof}
For a compactly supported Hamiltonian diffeomorphism $\psi$ of $\mathbb{R}^{2}$ that displaces the disk $\mathbb{D}$ one has $c_{-}(\phi) \geq -\gamma(\psi)$; see e.g. \cite{ginzburg-05}. The disk can be displaced by Hamiltonian diffeomorphisms with spectral norm arbitrary close to $1$ and so this proves the left inequality. The inequality on the right follows from $c_{+}(\phi) = -c_{-}(\phi^{-1})$. If $c_{+}(\phi) < 0$ then $c_{+}(\phi^{k})$ will be smaller than $-1$ for large enough $k$ by triangle inequality which is a contradiction; see item $(3)$ of Theorem \ref{T:Vit 1}. This argument can be found in \cite[Lemma 13]{Seyfaddini-15}. Hence $c_{+}(\phi) \geq 0$. Finally one uses $c_{-}(\phi) = -c_{+}(\phi^{-1})$ to finish the proof.
\end{proof}

In \cite{bramham2024spectral}, using compactly supported extensions for Hamiltonians explained in \S\ref{section:approx}, the authors extended the definition of Viterbo's spectral invariants 
to the isotopy class of arbitrary smooth Hamiltonians on $\R/\Z \times \D$. 
The construction is as follows. Let $\varphi : \mathbb{D} \to \mathbb{D}$ be an area-preserving diffeomorphism and
$\{\varphi_t\}_{t \in [0,1]}$ be an area-preserving isotopy from $id$ to $\varphi$. Let $H_t$ be a Hamiltonian function generating the isotopy $\varphi_t$ such that 
$H_t|_{\partial \mathbb{D}} = 0$. We choose a sequence of Hamiltonians $H_n$ compactly supported in a slightly larger disk, which converge to $H$ in the $C^0$–topology as in \S\ref{section:approx}. The Viterbo spectral invariants can be extended by setting
\begin{equation}\label{eq:spec_inv_non_comp_supp_maps}
\sigma_{\pm}(\varphi_{t}) = \sigma_{\pm}(H) := \lim_{n \to \infty} c_{\pm}(H_n).
\end{equation}
The authors showed that these extended spectral invariants satisfy all the properties 
listed in Theorem~\ref{T:Vit 1}, except for the spectrality property, 
which may fail in this general setting. However, a weaker version of spectrality still holds.
\begin{theorem}\label{T:spec main-1 Intro}
Let $\varphi_t$ be a Hamiltonian isotopy on the disk $\D$ with area $1$, not necessarily compactly supported in $\D$. Then, the extended spectral invariants $\sigma_{\pm}(\varphi_{t})$ depend only on the homotopy class of $\varphi_{t}$ with fixed endpoints, and 
\begin{equation}
\sigma_{\pm}(\varphi_{t})\in\spec(\varphi_{t})\cup\{0\}\cup\{\pm 1\},
\end{equation}
where $\spec(\varphi_{t})$ is the set of actions of the one-periodic orbits inside $\D$. Moreover, if $\rho(\varphi_{t}|_{\partial\mathbb{D}}) \in [0,1]$ then $\sigma_+(\varphi_{t})\in\spec(\varphi_{t})$. 
\end{theorem} 
The following example will be useful later in the proof of the main result.
\begin{example}
For isotopy $R_{2\pi t}(r,\theta):=(r,\theta+2\pi t)$ of the two-disk we have
\begin{enumerate}
\item $\sigma_{+}(R_{2\pi t})=0, \sigma_{-}(R_{2\pi t})=-1.$
\item $\sigma_{-}(R_{-2\pi t})=0, \sigma_{+}(R_{-2\pi t})=1.$
\end{enumerate}
Indeed, we have only one spectral value, which is $1,$ and $\rho(R_{2\pi t}|_{\partial\mathbb{D}})= 1$.
\end{example}

\subsubsection{The boundary depth}\label{sec:bourcode_bound_depth}
\begin{figure}
  \centering
  \begin{tikzpicture}[scale = 1]
\draw[->] (0,0) -- (0,3.5);
\draw[->] (0,0) -- (6,0);
\draw (0,0) -- (6,3);

\draw[dashed] (5,0) -- (5,2.5);
\draw[dashed] (0,2.5) -- (5,2.5);
\node at (5,-0.2) {\tiny{$1$}};

\draw[dashed] (4,0) -- (4,2);
\draw[dashed] (0,2) -- (4,2);
\node at (4,-0.2) {\tiny{$1 - \delta$}};

\draw[dashed] (2,0) -- (2,1.3);
\draw[dashed] (0,1.3) -- (2,1.3);
\node at (2,-0.2) {\tiny{$1 - 2\delta$}};

\draw[line width=0.3mm, color = black] (4,2) -- (6,3);
\draw[line width=0.3mm, color = black] (2,1.3) to[out = 0, in = -155, looseness = 1] (4,2);
\draw[line width=0.3mm, color = black] (0,1.3) -- (2,1.3);
\node at (-0.2,2.5) {\tiny{$\theta$}};
\node at (-0.6,2) {\tiny{$\theta (1-\delta)$}};
\node at (-0.8,1.3) {\tiny{$\theta(1 - 3\delta/2)$}};
\node at (6.3,0) {\tiny{$|z|^{2}$}};
\end{tikzpicture}
\caption{The function $H_{\theta, \delta}$.}
\label{fig:hamiltonian}
\end{figure}
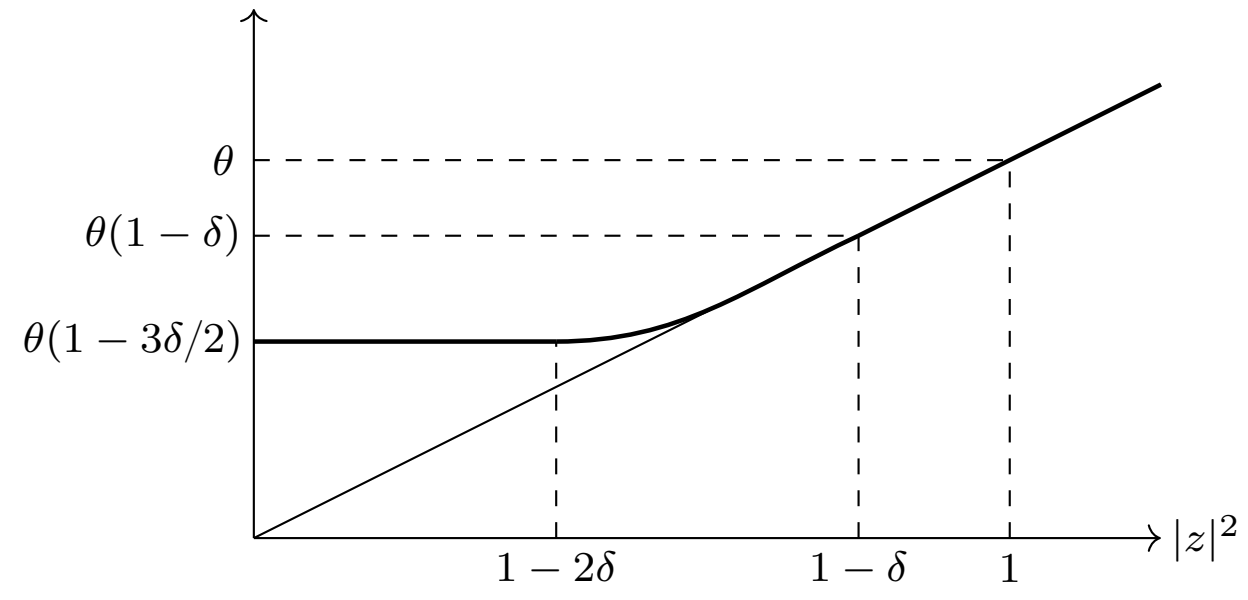

In this section we define the {\em boundary depth}, a notion that plays a central role in the proof of the main theorem. Consider the following Hamiltonian,
\begin{equation}\label{eq:the_hamiltonian}
H_{\theta,\delta}(z) = \theta(\mu_{\delta}(\vert z\vert^{2}-1)+ 1),
\end{equation}
where $\mu_{\delta}$ is a convex cut-off function so that:
\begin{enumerate}
\item $\mu_{\delta}(x)=x$ for $-\delta \leq x$,
\item $\mu_{\delta}(x) \equiv -3\delta/2$ for $x\leq -2\delta$,
\item $\mu_{\delta}'(x)>0$ for $-2\delta < x \leq -\delta$.
\end{enumerate}
It is important that $H_{\theta,\delta}(z) = \theta |z|^{2}$ near the boundary of the disk $\D$ where $\{1 - \delta \leq |z|^{2}\}$ and $H_{\theta,\delta} \equiv \theta(1 - 3\delta/2)$ is constant on $\{|z|^{2} \leq 1 - 2\delta\}$; see Figure \ref{fig:hamiltonian}.

Let $\phi^{1}_{\theta,\delta}$ be the time-one map of the flow generated by $H_{\theta,\delta}$; it is rotation by angle $2\pi \theta$ near the boundary of the disk $\D$ where $|z|^{2} \geq 1 - \delta$.

The boundary depth $\beta(\varphi)$ of a Hamiltonian diffeomorphism $\varphi$ of a closed symplectic manifold is studied in \cite{Usher-13} and is defined to be the maximal length of a finite bar in its associated barcode $\mathcal{B}(\varphi)$. To specify the coefficient field $\F$ used in the definition we use the notation $\mathcal{B}(\varphi, \F)$. 

The following is the definition of the boundary depth for compactly supported Hamiltonian diffeomorphisms of the the non-compact Liouville manifold $\C$ with the Liouville form $\lambda_{st} = \frac{r^{2}}{2\pi}d\theta$; see \cite{Frauenfelder-Schlenk-07} for Floer theory in the non-compact setting. The persistence module we use here is $V_{\phi,t} = \HF^{\leq t}(H)$ and the structure maps $V_{\phi,s} \rightarrow V_{\phi,t}$ for $s \leq t$ are maps induced by inclusions $\CF^{\leq s}(H)\rightarrow \CF^{\leq t}(H)$. 
\begin{definition}\label{def:bound_depth_liou}
Let $H$ be a Hamiltonian compactly supported inside $\mathbb{D}$ that generates the Hamiltonian diffeomorphism $\varphi$ as its time-one map. Suppose that its support lies in $\{|z|^{2} < 1 -2\delta\}$ for some small $\delta > 0$. Let $H_{\theta, \delta}$ be as in \eqref{eq:the_hamiltonian},  where $\theta = \epsilon$ for some $0 < \epsilon < 1$, generating $\phi^{1}_{\epsilon, \delta}$ as its time-one map. We define the boundary depth of $\varphi$ by
$$\beta(\varphi) := \beta(\varphi \circ \phi^{1}_{\epsilon, \delta}).$$
The boundary depth of $\varphi \circ \phi^{1}_{\epsilon, \delta}$ can be defined via a limiting process, since for every two non-degenerate \emph{compactly supported} perturbations $\psi_{1},\psi_{2}$ of $\varphi \circ \phi^{1}_{\epsilon, \delta}$ the argument for part $(\mathrm{iii})$ of \cite[Theorem 1.4]{Usher-13} goes through and,
$$|\beta(\psi_{1}, \mathbb{K}) - \beta(\psi_{2}, \mathbb{K})| \leq |\psi_{1}^{-1}\psi_{2}|_{\Hof};$$
note that $\psi_{1}$ and $\psi_{2}$ both have rotation angle $2\pi \epsilon$ near the boundary and so there are continuation maps in both directions
$$\mathrm{CF}(\psi_{1}) \overset{\mathfrak{c}_{12}}\longrightarrow \mathrm{CF}(\psi_{2}) \overset{\mathfrak{c}_{21}}\longrightarrow \mathrm{CF}(\psi_{1}),$$
defining a $|\psi_{1}^{-1}\psi_{2}|_{\Hof}$-quasiequivalence, see \cite[Definition 3.7]{Usher-13}, between the Floer complexes of $\psi_{1}$ and $\psi_{2}$. Here $|.|_{\Hof}$ is the Hofer norm.
\end{definition}


\begin{proposition}\label{clm:ind_bound_depth}
The definition of the boundary depth is independent of $\epsilon$ and $\delta$.
\end{proposition}
A careful proof of Proposition \ref{clm:ind_bound_depth} is written in Appendix \ref{app:bound_depth_well_defined}. 
\begin{lemma}\label{lem:cont_of_bound_depth}
For every compactly supported Hamiltonian diffeomorphism $\varphi$ of the disk $\D$ one has $\beta(\varphi, \mathbb{K}) \leq 1,$ where $\mathbb{K}$ is any ground field.
\end{lemma}
\begin{proof}
This follows from the inequality $\beta \leq \gamma$ proved in \cite[Theorem 1.10]{feng-zhang} for Liouville domains and item $(5)$ of Theorem \ref{T:Vit 1}, see also \cite[Theorem A]{kislev2022bounds}, and the displacement energy bound on the spectral norm \cite{Schwarz-00, ginzburg-05} providing the bound $\gamma(\varphi)\leq 2$.
\end{proof}

\begin{proposition}\label{lem:finite_boundary_depth}
For every Hamiltonian diffeomorphism $\varphi$ of the disk that is linear at infinity, not necessarily with small slope, one has $\beta(\varphi, \mathbb{K}) \leq 1.$
\end{proposition}

\begin{proof} By Lemma \ref{lem:cont_of_bound_depth} it is enough to show that for every Hamiltonian diffeomorphism $\varphi$ of the disk that is linear at infinity with slope $\theta$ there is a Hamiltonian diffeomorphism $\varphi'$ linear at infinity with slope $0 \leq \theta' < 1$ so that
$$\beta(\varphi, \mathbb{K}) = \beta(\varphi', \mathbb{K}).$$ 
Write $\varphi = \phi \circ \phi^{1}_{\theta, \delta}$ for some small $0 < \delta< 1$ where $\phi^{1}_{\theta, \delta}$ is generated by $H_{\theta, \delta}$, see \eqref{eq:the_hamiltonian}, and coincides with the rotation map by angle $2\pi\theta$ on $\{1-\delta \leq |z|^{2}\}$ and $\phi$ is compactly supported in $\mathrm{int}(\mathbb{D})$. Consider the Hamiltonian $K(z) = |z|^{2}$ which generates a loop of Hamiltonian diffeomorphisms. There is a naturality transformation associated with a suitable integer multiple of $K$:
$$\mathfrak{n}: \mathrm{CF}(\phi \circ \phi^{1}_{\theta, \delta}) \rightarrow \mathrm{CF}(\phi \circ \phi^{1}_{\epsilon, \delta}),$$
where $0 \leq \epsilon < 1$ and $\epsilon = \theta\in \R/\Z$; there are some choices of almost complex structures which we avoid mentioning here as they are irrelevant; see \cite[\S5]{shelukhin-zhao-22} for more details. In other words, the naturality transformation reduces the rotation angle of $\varphi$ near the boundary from $2\pi \theta$, which can be potentially large, to a small rotation angle $0 \leq 2\pi \epsilon < 2\pi$. It is well known that the (ungraded) Floer complexes are identical and the transformation $\mathfrak{n}$ is just the identity map on the underlying vector spaces and preserves actions; so does its inverse $\mathfrak{n}^{-1}$; see for instance \cite[\S2.8]{alizadeh-atallah-cant-ellipsoid-24}. Hence for $\varphi' = \phi \circ \phi^{1}_{\epsilon, \delta}$ we have  
$$\beta(\varphi, \mathbb{K}) = \beta(\varphi',\mathbb{K})$$
and $\varphi'$ is linear at infinity with angle $0 \leq 2\pi \epsilon < 2\pi$.
\end{proof}

\subsubsection{Smith inequality for Floer homology}\label{sec:smith_ineq_ham_floer_cohom}
Here we recall a key Smith-type inequality for Floer homology proved in \cite[Theorem C]{shelukhin-zhao-22} and some of its useful consequences.

\begin{theorem}\label{thm:smith_ineq} Let $\varphi$ be a Hamiltonian diffeomorphism of $\mathbb{D}$ that is a rotation near the boundary. Let $p$ be a prime number and $I$ be an interval with end points of $I$ and $pI$ being disjoint from the spectrum of $\varphi$. Then, the following inequality holds:
\begin{equation}\label{eq:smith-ineq}
\dim \mathrm{HF}(\varphi, \mathbb{F}_{p})^{I} \leq \dim \mathrm{HF}(\varphi^{p}, \mathbb{F}_{p})^{pI},
\end{equation}
where the dimensions are over the field $\mathbb{F}_{p}$ with $p$ elements.
\end{theorem}
The following corollary will be used later in the proof of the main theorem.
\begin{corollary}\label{rem:bars_covered_by_finite_bars}
Let $m(t,\varphi)$ be the number of finite bars $I$ in the barcode of $\varphi$ with $t\in I$. Then the following holds for all primes $p$:
$$m(t, \varphi) \leq m(p t, \varphi^{p}).$$
Moreover, if $I \in \mathcal{B}(\varphi, \mathbb{F}_{p})$ is a finite bar then the interval $pI$ will be covered by a collection of finite bars from the barcode $\mathcal{B}(\varphi^{p}, \mathbb{F}_{p})$ of $\varphi^{p}$. 
\end{corollary}

\begin{proof}
Let $\underline{n}(t,\varphi)$ the number of infinite bars starting before $t$ and $\overline{n}(t,\varphi)$ be the number of infinite bars starting after $t$ in the barcode $\mathcal{B}(\varphi, \mathbb{F}_{p})$ of $\varphi$. Consider the Smith-type inequality \eqref{eq:smith-ineq} for the intervals $(-\infty, t)$ and $(t, \infty)$:
$$\eqref{eq:smith-ineq} \text{ for } (-\infty, t) \implies m(t,\varphi) + \underline{n}(t, \varphi) \leq m(pt,\varphi^{p}) + \underline{n}(pt, \varphi^{p}),$$
$$\eqref{eq:smith-ineq} \text{ for } (t,+\infty) \implies m(t,\varphi) + \overline{n}(t, \varphi) \leq m(pt,\varphi^{p}) + \overline{n}(pt, \varphi^{p}).$$
Summing up the two inequalities it is enough to prove the following equality,
$$\underline{n}(t, \varphi) + \overline{n}(t, \varphi) = \underline{n}(pt, \varphi^{p}) + \overline{n}(pt, \varphi^{p}).$$
It is well-known that the disk $\mathbb{D}$ admits a loop of Hamiltonian diffeomorphisms generated by $K(z) = |z|^{2}$. This loop and its iterations, induce filtered isomorphisms, so called naturality transformations, between Floer cohomologies $\mathrm{HF}(\varphi^{p}, \mathbb{F}_{p})$, $p$ any prime number, and Floer cohomologies of Hamiltonian diffeomorphisms with slope in $(0,1)$. In particular the barcodes coincide. The latter are isomorphic with cohomology of the disk and have dimensions equal to one. Hence, their barcode, and consequently barcodes $\mathcal{B}(\varphi^{p}, \mathbb{F}_{p})$ have only a single infinite bar. Therefore,
$$\underline{n}(t, \varphi) + \overline{n}(t, \varphi) = 1 = \underline{n}(pt, \varphi^{p}) + \overline{n}(pt, \varphi^{p}),$$
which proves the inequality. To prove the second part, take an interval $I \in \mathcal{B}(\varphi, \mathbb{F}_{p})$. Let $t \in pI$. Since $I$ contains $t/p$ it follows that, 
$$1 \leq m(t/p,\varphi) \leq m(t, \varphi^{p}),$$
and hence there is a finite bar in $\mathcal{B}(\varphi^{p}, \mathbb{F}_{p})$ containing $t$.
\end{proof}

\subsection{Franks-Handel theorem on zero entropy maps of \texorpdfstring{$S^{2}$}{}}
In this section, we introduce the classification of zero-entropy maps given by Franks and Handel 
in \cite{franks2013entropy} for area-preserving maps of $S^{2}$ that are smooth on the complement of their fixed points.  
Later, this was improved by Le Calvez and Tal in \cite{le2022topological} for area-preserving homeomorphisms without horseshoes, 
showing that there does not exist an invariant subset $\Lambda$ as in Theorem~\ref{horseshoe}.  
This result will be crucially used in the highly degenerate case where Floer theory is difficult to apply. In the following theorem $h_{top}(f)$ is the topological entropy of $f$, see Definition \ref{def:top_entropy}, $\Omega(f)$ is the set of non-wandering points of $f$, i.e., 
$$\Omega(f) = \{x\in \D: \forall U_{x}\ \exists k\ \text{ s.t. } f^{k}(U_{x}) \cap U_{x} \neq \emptyset\},$$
where $U_{x}$ is a neighborhood of $x$, the $\omega$-limit set of a point $x\in \D$ is the set of limit points of the sequence $\{f^{k}(x)\}_{k\geq 1}$ and is denoted by $\omega(x)$, and finally the $\alpha$-limit set of $x$ is defined to be the set of converging points of the sequence $\{f^{-k}(x)\}_{k\geq 1}$ and is denoted by $\alpha(x)$. The word ``\emph{frontier}" of an open annulus $A_{\alpha}$ in the following theorem will be referring to $\partial\bar{A}_{\alpha}\backslash \mathrm{Fix}(f)$ where the latter is the fixed point set of $f$. The following is proved in \cite[Theorem 1.2, Proposition 5.1]{franks2013entropy}.
\begin{theorem}\label{thm:franks_handel_theorem}
Let $f : S^{2} \to S^{2}$ be an area-preserving homeomorphism with zero topological entropy, infinite order and at least three periodic points, and is smooth away from a finite set of singular fixed points. Then, there exists a family of pairwise disjoint invariant open sets $\mathcal{A}(f) = \{A_{\alpha}\}_{\alpha}$ whose union is dense in the complement of $\mathrm{Fix}(f)$, and, moreover,
\begin{enumerate}
    \item the sets $A_{\alpha}$ are the maximal fixed-point-free invariant open regular annuli; ``regular" means that each annulus is the interior of its closure in $S^{2}\backslash \Fix(f)$; \cite[Lemma 14.5]{franks2013entropy}.
    \item the $\alpha$-limit set of every point $z \notin \cup A_{\alpha}$ is contained in a single connected component of $\mathrm{Fix}(f)$; the same holds for the $\omega$-limit set of $z$.
    \item Let $C$ be a connected component of the frontier of $A_{\alpha}$. Then, the connected components of $\mathrm{Fix}(f)$ that contain $\alpha(z)$ and $\omega(z)$ are independent of $z \in C$.
\end{enumerate}
\end{theorem}
The elements of $\mathcal{A}(f)$ are the components of the set of weakly free disk recurrent points of $f$, see \cite[Definition 1.1]{franks2013entropy}, which is an open dense subset of $S^{2}\backslash \Fix(f)$ containing all the bi-recurrent points.
\begin{remark}
We can apply this result to an area-preserving map of the disk because by collapsing the boundary of the disk to a point, we obtain an area-preserving homeomorphism of $S^{2}$.  The collapsed point becomes a fixed point of this map, and the map is smooth on the complement of this fixed point.
\end{remark}
We know that if $f$ has zero topological entropy then so do all of its iterations. Hence, for every $q$ we can apply Theorem \ref{thm:franks_handel_theorem} to $f^{q}$ and find a collection of disjoint invariant open annuli that are dense in the complement of $\mathrm{Fix}(f^{q})$ and satisfy the extra properties of Theorem \ref{thm:franks_handel_theorem}; we denote this collection by $\mathcal{A}(q)$ if there is no confusion on the choice of the map $f$. Franks and Handel \cite{franks2013entropy} showed that these annuli satisfy the following additional properties.
\begin{theorem}\label{thm:regopen}
The annuli $\mathcal{A}(q)$ given by Theorem \ref{thm:franks_handel_theorem} for $f^{q}$ satisfy the following:
\begin{enumerate}
\item each $A \in \mathcal{A}(q)$ is the interior of its closure.
\item The map $f$ permutes the elements of $\mathcal{A}(q)$.
\item Each $A\in \mathcal{A}(q)$ is a subset of some $A' \in \mathcal{A}(1)$.
\end{enumerate}
\end{theorem}
Let $A = S^{1} \times \mathbb{R}$  be the open annulus and $f$ be a homeomorphism of $A$ isotopic to the identity.  A point $z$ is called an \emph{escaping point} of $f$ if the sequence 
$\{f^{n}(z)\}_{n\geq 1}$ diverges to infinity, i.e., its 
$\omega$-limit set $\omega(z)$ is empty. We denote by $ \mathrm{ne}^{+}(f)$ the complement of the set of escaping points; it is a closed subset of $A$. We also define $\mathrm{ne}^{-}(f) := \mathrm{ne}^{+}(f^{-1})$; it is consists of points whose $\alpha$-limit set is nonempty. Finally, we set,
$$\mathrm{ne}(f) =  \mathrm{ne}^{+}(f) \cup  \mathrm{ne}^{-}(f).$$
The following version of the continuity of the rotation number was proven by Le Calvez and Tal \cite{le2022topological}.
\begin{theorem}\label{thm:rotcont}
Let $f$ be a homeomorphism of an annulus $A$ isotopic to $id$. Let $\widehat{f}$ be its lift to $\mathbb{R}^2$. Moreover, suppose that $f$ has no topological horseshoe. Then:
\begin{enumerate}
    \item each point $z \in \mathrm{ne}^{+}(f)$ has a well-defined rotation number $\mathrm{rot}_{\widehat{f}}(z)$.
    \item For every $z, z_0 \in \mathrm{ne}^{+}(f)$ such that $z_0 \in \omega(z)$, we have $\mathrm{rot}_{\widehat{f}}(z_0) = \mathrm{rot}_{\widehat{f}}(z).$
    \item If $z \in \mathrm{ne}^{+}(f) \cap \mathrm{ne}^{-}(f)$ is non-wandering, then $\mathrm{rot}_{\widehat{f}^{-1}}(z) = - \mathrm{rot}_{\widehat{f}}(z).$
    \item The map $\mathrm{rot}_{\widehat{f}^{\pm}} : \Omega(f) \cap \mathrm{ne}(f) \to \mathbb{R}$ defined below is continuous,
    $$\mathrm{rot}_{\widehat{f}^{\pm}}(z) =
        \begin{cases}
            \mathrm{rot}_{\widehat{f}}(z) & \text{if } z \in \Omega(f) \cap \mathrm{ne}^{+}(f),\\[2mm]
            -\mathrm{rot}_{\widehat{f}^{-1}}(z) & \text{if } z \in \Omega(f) \cap \mathrm{ne}^{-}(f).
        \end{cases}
    $$
\end{enumerate}
\end{theorem}
See the Example after \cite[Theorem A]{le2022topological} for an example of an area preserving homeomorphism $f$ isotopic to the identity for which the rotation number functions $\mathrm{rot}_{\widehat{f}^{\pm 1}}$ are not continuous on their domains of definition $\mathrm{ne}^{\pm}(f)$ and the equality $\mathrm{rot}_{\widehat{f}} = -\mathrm{rot}_{\widehat{f}^{-1}}$ holds only on non-wandering points of $f$.

\begin{remark}\label{rem:rest_are_isotop_to_id}
Throughout \S\ref{sec:maps_with_prop_p_main_thm} we often work with the restriction of a symplectic diffeomorphism $\phi$ of a disk to an invariant open topological annulus $A$. To freely apply results of \cite{franks1988generalizations, franks2013entropy} to $\phi|_{A}$ we need to confirm that the prime extension $(\phi|_{A})_{c}$ of $\phi|_{A}$ to the prime end compactification $A_{c}$ of $A$ is isotopic to the identity. By an immediate topological reason we know that $\phi$ preserves the boundary components of $A$ and so $(\phi|_{A})_{c}$ preserves the boundary components of $A_{c}$. Moreover, $(\phi|_{A})_{c}$ is an orientation preserving homeomorphism of the closed annulus $A_{c}$. So $(\phi|_{A})_{c}$ induces identity on the homology group $\H(A_{c}, \partial A_{c}, \Z)$ and it is \emph{homotopic} to the identity, which is immediate by linear interpolation in the universal cover. Now it follows from a theorem of Baer \cite{MR1581099, MR1581156} and Epstein \cite{MR214087}, see also \cite[Theorem 1.12]{MR2850125}, that $(\phi|_{A})_{c}$ is isotopic to the identity via homeomorphisms.
\end{remark}

\section{Proof of the main theorem}
We prove Theorem~\ref{mainthm:meanacion} by dividing the argument into several cases. First, we prove the statement for compactly supported maps. Building on this, we extend the proof to maps with rational rotation number on the boundary. Finally, we consider maps with irrational rotation number on the boundary; this case will be further divided into three sub-cases. In this analysis, maps with positive and zero topological entropy play a crucial role.

\subsection{Compactly supported}\label{sec:proof_comp_supp}
We first prove Theorem \ref{mainthm:meanacion} for compactly supported Hamiltonian diffeomorphisms.
\begin{theorem}\label{thm:compmean}
Let $\phi \neq id$ be a compactly supported area-preserving diffeomorphism of $\D$. Then, for an interval $[a,b] \subset [c_{-}(\phi), 0)$, every prime number $p,$ and positive integer $k$ with
$$p^{k} > \max\{\frac{1}{\vert b\vert}, \frac{1}{\vert b - a\vert}\},$$
there exist an $p^{k}$-periodic orbit with mean action in the interval $[a,b]$. A similar result holds for the interval $(0, c_{+}(\phi)] = -[c_{-}(\phi^{-1}), 0)$. Moreover, the closure of the mean action spectrum of $\phi$ contains all the finite bars in the $\F_{p}$-barcodes of $\phi^{p^{k}}$ scaled by $p^{-k}$ for all prime numbers $p$ and positive integers $k$.
\end{theorem}
\begin{proof}
The spectral invariants $c_{\pm}(\phi)$, and the action spectrum $\mathrm{Spec}(\phi)$, the critical values of $\mathcal{A}_{H}$ where $H\in C^{\infty}_{c}(\R/\Z\times \D)$ generates $\phi$ as its time-one map, are continuous in $H\in C^{\infty}_{c}(\R/\Z\times\D)$ where the latter is equipped with $C^{\infty}$-topology. Given a fixed number $\ell\geq 1$ the map 
\begin{equation}\label{eq:cont_spectra}
H\mapsto (c_{\pm}(\phi), c_{\pm}(\phi^{2}), \dots, c_{\pm}(\phi^{\ell}), \mathrm{Spec}(H),\frac{1}{2}\mathrm{Spec}(H^{(2)}), \dots, \frac{1}{\ell}\mathrm{Spec}(H^{(\ell)})),
\end{equation}
is also continuous in $H$ where the spectra $\mathrm{Spec}(H^{(i)})$ are considered as subsets of $\R$ equipped with the Hausdorff distance. So it is enough to prove the theorem for a strongly non-degenerate perturbation $\varphi$ of $\phi$ that is an irrational rotation near the boundary $\partial\D$ with small positive irrational angle, i.e., we fix an isotopy $\varphi_{t}$ with $\varphi_{1} = \varphi$ generated by a Hamiltonian radially-linear near $\partial\D$ with small positive irrational slope; note that the choice of the Hamiltonian does not matter as any two choices will give \emph{filtered} isomorphic Floer homologies.

Since $\varphi \neq id$ and $\gamma: \Ham(\D)\rightarrow \R$ is a norm \cite{Viterbo-99, Schwarz2000, Oh2005} we have 
$$\gamma(\varphi) = c_{+}(\varphi) + c_{-}(\varphi) > 0.$$
Recall that $c_{-}(\varphi)\leq 0\leq c_{+}(\varphi)$. Without loss of generality assume that $c_{-}(\varphi) < 0$, otherwise one has $c_{+}(\varphi) > 0$ and one can instead work with $\varphi^{-1}$ which satisfies $c_{-}(\varphi^{-1}) = -c_{+}(\varphi) < 0$. Note that an $p^{k}$-periodic orbit of $\varphi^{-1}$ with mean action inside $[c_{-}(\varphi^{-1}), 0)$ corresponds to an $p^{k}$-periodic point of $\varphi$ with mean action inside $(0, c_{+}(\varphi)]$. Hence, assume $c_{-}(\varphi) < 0$, and consider $[a,b] \subset [c_{-}(\varphi), 0)$. Let $p$ be a prime number and $k$ a positive integer satisfying the condition of the statement. It is clear that,
$$\dim \HF^{[-1 - \epsilon, b]}(\varphi, \F_{p}) \geq 1,$$
for some small $\epsilon > 0$; recall that $-1 \leq c_{-}(\phi)$ by Lemma \ref{lem:pos_spec_inv}, and so $-1 - \epsilon < c_{-}(\varphi)$ must hold taking $\varphi$ small enough perturbation compared to $\epsilon$. It follows from the Smith inequality for Floer homology, Theorem \ref{thm:smith_ineq}, that,
$$\dim \HF^{[p^{k}(-1 - \epsilon), p^{k}b]}(\varphi^{p^{k}}, \F_{p}) \geq 1.$$
By the assumption $p^{k}\vert b\vert > 1$ and the fact that $-1 \leq c_-(\phi^{p^k})\leq 0$, see Lemma \ref{lem:pos_spec_inv}, we conclude the existence of a finite bar $[x,y]$ in the $\F_{p}$-barcode of $\varphi^{p^{k}}$ with $x < p^{k}b < y$. It follows from the assumption $p^{k}\vert b - a\vert > 1$ and Lemma \ref{lem:finite_boundary_depth} that $p^{k}a < x$. By spectrality of the ends of the bars \cite{Usher-13} we know that $x$ represents the action of an $p^{k}$-periodic orbit of $\varphi$, hence, the interval $[a,b]$ contains the mean action of an $p^{k}$-periodic orbit of $\varphi$.
\begin{figure}[h]
\begin{tikzpicture}[scale = 1]
\begin{scope}
\draw (0, 0) -- (4, 0);

\draw (0.5, 0.15) -- (0.5, -0.15) node[below=2pt] {\tiny{$-1$}};
\draw (1.2, 0.15) -- (1.2, -0.15) node[below=2pt] {\tiny{$c_{-}$}};
\draw (2.8, 0.15) -- (2.8, -0.15) node[below=2pt] {\tiny{$0$}};

\begin{scope}
    \clip (1.6, -0.8) rectangle (2.4, 0.8);
    \foreach \x in {1.4, 1.5, ..., 2.5} {
        \draw (\x, -0.1) -- (\x+0.25, 0.1);
    }
\end{scope}
\draw (1.6, -0.12) rectangle (2.4, 0.12);

\fill (1.2, 0.7) circle (1.5pt);
\draw (1.2, 0.7) -- (3.5, 0.7);
\end{scope}

\begin{scope}[xshift =1cm]
\draw[->] (3.7, 0.05) -- (5.5, 0.05);
\node at (4.8, 0.45) {$p^{k}$-th};
\node at (4.8, -0.4) {power};
\end{scope}

\begin{scope}[xshift = 2cm]
\draw (5.5, 0) -- (13.5, 0);

\begin{scope}
    \clip (6.3, -0.8) rectangle (8.8, 0.8);
    \foreach \x in {6.0, 6.1, ..., 8.9} {
        \draw (\x, -0.12) -- (\x+0.3, 0.12);
    }
\end{scope}
\draw (6.3, -0.15) rectangle (8.8, 0.15);

\draw (9.8, 0.15) -- (9.8, -0.15) node[below=2pt] {\tiny{$-1$}};
\draw (10.6, 0.15) -- (10.6, -0.15) node[below=2pt] {\tiny{$c_{-}(\varphi^{p^{k}})$}};
\draw (11.8, 0.15) -- (11.8, -0.15) node[below=2pt] {\tiny{$0$}};

\fill (5.8, 1.2) circle (1.5pt);
\draw[dashed] (5.8, 1.2) -- (13.5, 1.2);
\draw[dashed] (5.8, 0) -- (5.8, 1.2);
\node at (5.8, -0.2) {\tiny{$p^{k}c_{-}$}};
\fill (6, 0.6) circle (1pt);
\draw[ForestGreen] (6, 0.6) -- (6.6, 0.6);
\fill (6.6, 0.6) circle (1pt);

\fill (6.4, 0.4) circle (1pt);
\draw[ForestGreen] (6.4, 0.4) -- (7.3, 0.4);
\fill (7.3, 0.4) circle (1pt);

\fill (7.1, 0.6) circle (1pt);
\draw[ForestGreen] (7.1, 0.6) -- (7.8, 0.6);
\fill (7.8, 0.6) circle (1pt);

\fill (7.5, 0.4) circle (1pt);
\draw[ForestGreen] (7.5, 0.4) -- (8.7, 0.4);
\fill (8.7, 0.4) circle (1pt);

\fill (8.3, 0.6) circle (1pt);
\draw[ForestGreen] (8.3, 0.6) -- (9.0, 0.6);
\fill (9.0, 0.6) circle (1pt);

\fill (5.8, -0.6) circle (1.5pt);
\draw[red] (5.8, -0.6) -- (9.0, -0.6);
\fill (9.0, -0.6) circle (1.5pt);

\fill (10.6, 0.8) circle (1.5pt);
\draw[ForestGreen] (10.6, 0.8) -- (13.5, 0.8);
\end{scope}
\end{tikzpicture}
\caption{The barcode of $\varphi$ on the left has only one infinite bar starting by definition at $c_{-} = c_{-}(\varphi)$. The box on left is $[a,b]$. The barcode of $\varphi^{p^{k}}$ on the right has only one infinite bar starting at $c_{-}(\varphi^{p^{k}})\in [-1,0]$. It has dimension at least one at every point of the large box $[p^{k}a, p^{k}b]$ on the right but long bars such as the red one are not allowed by the estimate $\beta(\varphi^{p^{k}})\leq 1$.}
\end{figure}
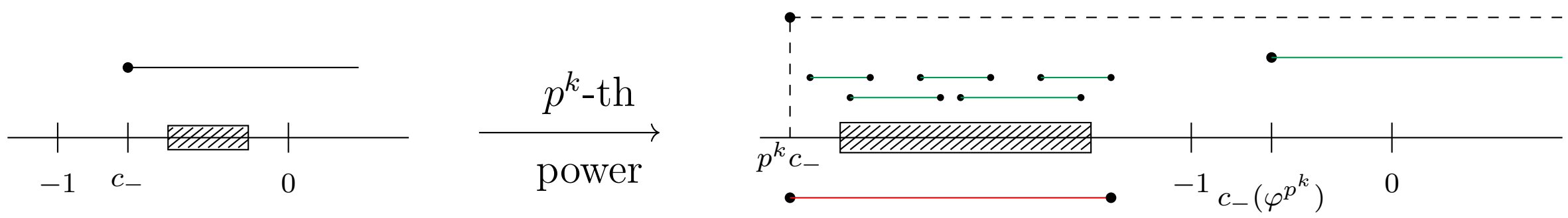
To prove the second part of the theorem assume that $[a,b]$ is a finite bar in the $\F_{p}$-barcode of $\varphi$; the case of $\varphi^{p^{k}}$ stated in the theorem is similar because an $\ell$-periodic point of $\varphi^{p^{k}}$ with mean action equal to $x\in [a,b]$ is an $\ell p^{k}$-periodic point of $\varphi$ with mean action $p^{-k}x\in [p^{-k}a, p^{-k}b]$. By Corollary \ref{rem:bars_covered_by_finite_bars} the interval $[p^{k}a, p^{k}b]$ is covered by finite bars $\cup_{k}J_{k}$ from the $\F_{p}$-barcode of $\varphi^{p^{k}}$. The length of each bar $J_{k}$ is bounded from above by $1$; see Lemma \ref{lem:finite_boundary_depth}. Since the ends of each bar $J_{k}$ represent actions of periodic orbits of $\varphi$, by increasing $k$ we observe that the mean action spectrum of $\varphi$ is dense in the interval $[a,b]$.
\end{proof}

\subsection{Rational rotation number}\label{sec:case_rational_rot_num}
We prove Theorem \ref{mainthm:meanacion} in the case where the rotation number on the boundary is rational. By passing to an iteration and supposing that the map is not periodic, we may and will assume that the map is non-identity with integer rotation number on the boundary.
\begin{theorem}\label{thm:ratcase}
Let $\phi\neq id$ be an area-preserving diffeomorphism of $\D$ with integer rotation number on the boundary. There exist numbers $-1\leq c_{-}\leq 0\leq c_{+}\leq 1$ so that $c_{+} - c_{-} > 0$ and the following holds: for any interval $[a,b] \subset [c_{-}, 0)$, prime number $p$ and positive integer $k$ so that,
$$p^{k} > \max\{\frac{1}{\vert b\vert}, \frac{1}{\vert b - a\vert}\},$$
there is an $p^{k}$-periodic orbit of $\phi$ with mean action in the interval $[a,b]$. Similar result holds for the intervals $[a,b] \subset (0, c_{+}]$ and $p^{k} > \max\{1/a, 1/(b-a)\}$.
\end{theorem}

\begin{proof}
In the following, by mean action of a periodic point $x$ of period $k$, for a Hamiltonian $G$, we mean $\mathcal{A}_{G^{(k)}}(x)/k$; see \eqref{eq:ham_act_funct} and compare it with $a_{(\phi_{G}^{t})^{k},\lambda}(x)/k$ as in Remark \ref{rem:compare_two_action_defs}. By a composition with a suitable full rotation isotopy of the disk, there exists a Hamiltonian isotopy $\phi_{t}$ with $\phi_{1} = \phi$ whose winding number on the boundary is zero. Let $H$ be a Hamiltonian generating $\phi_{t}$ so that $H \equiv 0$ on $\partial \D$. Let $H_{n}$ be the compactly supported Hamiltonians in a slightly larger disk $\D_{1 + 1/n}$ of radius $1 + 1/n$ approximating $H$ as in \S\ref{section:approx} and denote their Hamiltonian flows by $\phi_{n,t}$. Define,
$$c_{\pm} = \lim_{n\to \infty}c_{\pm}(H_{n}).$$
It follows immediately from the construction of $c_{\pm}$, item $(4)$ of Theorem \ref{T:Vit 1} and Lemma \ref{lem:pos_spec_inv} that $-1 \leq c_{-} \leq 0 \leq c_{+} \leq 1$ and $c_{+} - c_{-} > 0$. Assume that $c_{-} < 0$ and let $[a,b] \subset [c_{-}, 0)$ be some interval. Let $p$ and $k$ be as in the statement. By the result of \S\ref{sec:proof_comp_supp} we know that there is an $p^{k}$-periodic point $z_{n}\in \D_{1 + 1/n}$ of $\phi_{n,1}$ with mean action inside $[a,b]$ for every large enough $n$. 

The mean actions of $p^{k}$-periodic points of $\phi_{n,1}$ in $\D_{1 + 1/n}\backslash\D$ get arbitrary close to zero as $n\to \infty$; recall that the winding number of $\phi_{t}$, and so that of $p^{k}$-th iteration of $\phi_{t}$, are zero on $\partial \D$, so by Remark \ref{ext:smallact}, we obtain
$$a_{\phi_{n,1}^{p^{k}}, \lambda}(x)/p^{k} \to 0, \text{ as } n \to \infty,$$ 
for every $p^{k}$-periodic point $x\in \D_{1 + 1/n}\backslash\D$ of $\phi_{n,1}.$ Consequently by Remark \ref{rem:compare_two_action_defs} the mean actions,
$$\mathcal{A}_{H_{n}^{(p^{k})}}(x)/p^{k} = - a_{\phi_{n,1}^{p^{k}}, \lambda}(x)/p^{k},$$ 
converge to zero. We conclude that $z_{n}\in \D$ is a periodic point of $\phi$ for large enough $n$. The case of $c_{+} > 0$ and intervals $[a,b] \subset (0, c_{+}]$ is similar.
\end{proof}
\subsection{Irrational rotation number}\label{sec:case_irrational_rot_num}
We first shortly explain how to reduce this case to the case of rigid rotation near the boundary. Let $\psi: \D \rightarrow \D$ be a smooth area-preserving diffeomorphism with irrational rotation number on the boundary. Fix an open interval $(a,b)$ and an integer $k$. Recall that:
\begin{enumerate}
\item For every $n > k$, by Theorem \ref{thm:rotonboundary}, there is a smooth approximation $\psi_{n}$ of $\psi$ whose set of $k$-periodic points coincides with that of $\psi$ and their mean action functions are in distance $1/n.$
\item By Theorem \ref{thm:rotonboundary} $(3)$, there is a smooth area preserving diffeomorphism $\nu_{n}$ such that $\nu_{n}\psi_{n}\nu_{n}^{-1}$ is rigid rotation $\theta \mapsto \theta + 2\pi\theta_{0}$ on the boundary.
\item By Theorem \ref{thm:extensionrot}, there is a smooth extension $\phi_{n}$ of $\nu_{n}\psi_{n}\nu_{n}^{-1}$ to the disk $\D_{1 + 1/n}$ that is rigid rotation $\theta \mapsto \theta + 2\pi\theta_{0}$ near the boundary, with no new $k$-periodic points and whose mean action function is that of $\nu_{n}\psi_{n}\nu_{n}^{-1}$ shifted by a constant $a_{n}$ with $\vert a_{n}\vert \leq 3\theta_{0}/n$. Moreover, the mean action of periodic points in $\D_{1+1/n}\backslash\D$ are in $1/n$ distance from $\mathrm{area}(\D_{1 + 1/n})\theta_{0}$.
\end{enumerate}
Suppose that there is an $k$-periodic point $p\in \D_{1 + 1/n}$ of $\phi_{n}$ with mean action inside the interval $(a + a_{n},b + a_{n})$. By item $(3)$ above we must have $p\in \D$ being an $k$-periodic point of $\nu_{n}\psi_{n}\nu_{n}^{-1}$ with mean action inside $(a,b)$. Since the \emph{mean action} is preserved under conjugation by an area preserving diffeomorphisms, $\nu_{n}^{-1}(p)$ is an $k$-periodic point of $\psi_{n}$ with mean action in $(a,b)$. Since, for every fixed number $\ell$, the mean action spectrum $\cup_{1\leq i\leq \ell} \frac{1}{i}\mathrm{Spec}(H^{(i)})$ is continuous in $H$, the compactly supported Hamiltonian generating $\phi$, in the sense of \eqref{eq:cont_spectra}, it follows from item $(1)$ that if $n$ is sufficiently large then there is an $k$-periodic point of $\psi$ with mean action inside the interval $(a,b)$.

Hence, since the Smith inequality provides periodic points with mean action in a given interval $(a,b)$ whose period only depends on the interval itself, without loss of generality we can and will assume that the map is a rigid irrational rotation near the boundary by the discussion above.

We divide this case into three sub-cases. First, we assume:
\begin{enumerate}
\item $\phi: \D \rightarrow \D$ has positive topological entropy.
\end{enumerate}
In this case $\phi$ has many hyperbolic periodic points, up to replacing $\phi$ with its iteration we might assume that $\phi$ has many hyperbolic fixed points; see \S\ref{sec:positive_top_entrop}. Hence, there is at least one finite bar in the barcode of $\phi$, which we denote by $J$, such that at least one of its ends corresponds to a hyperbolic fixed point $p$ of $\phi$. 

It was shown in \cite[Theorem 3.1]{ginzburg-gurel-hyperbolic} that any Floer cylinder asymptotic to the hyperbolic one-periodic orbit corresponding to $p$, or to any of its iterations, has energy at least $c = c(p) > 0$ that depends only on the fixed point $p$, and not on its iterations; see also \cite[Theorem 6.1]{ginzburg-gurel-invent}. The idea is that, first, one reduces to the case where $p$ is a constant orbit \cite[\S9]{salamon-zehnder-92}, then, if there is a sequence of Floer cylinders, for possibly different iterations of the Hamiltonian, asymptotic to $p$ with energy converging to zero, via a Gromov-type convergence \cite{fish-11-gromov-comp}, one can find an integral curve of the Hamiltonian flow that stays inside a fixed small neighborhood of $p$ for a long enough time, and this is in contradiction with hyperbolicity of $p$.  

Continuing with the proof of case $(1)$, since the length of $J$ is at least energy of some Floer cylinder asymptotic to the orbit of $p$ we conclude that length of $J$ is non-zero and at least $c > 0$. Finally, similar analysis as in the proof of the second part of Theorem \ref{thm:compmean} shows that the closure of the mean action spectrum of $\phi$ contains the interval $J$.

Recall that $\rho(\phi_{t})$ denotes rotation number of $\phi_{t}$ on $\partial \D$. For the second case we assume:
\begin{enumerate}[resume]
\item $\phi: \D \rightarrow \D$ has zero topological entropy. Moreover, for every iteration $\phi^{n}$ and an isotopy $\phi_{n,t}$ with time-one map $\phi^{n}$ satisfying $\vert \rho(\phi_{n,t})\vert \leq 1$ we have,
\begin{enumerate}
\item $\sigma_{+}(\phi_{n,t})= \max\{-\rho(\phi_{n,t}), 0\}$.
\item $\sigma_{-}(\phi_{n,t}) = \min\{-\rho(\phi_{n,t}), 0\}$.
\end{enumerate}
\end{enumerate}
Here $\sigma_{\pm}$ are the spectral invariants defined for non-compactly supported area-preserving diffeomorphisms via compactly supported extensions; see \eqref{eq:spec_inv_non_comp_supp_maps}. This case is treated in Theorem \ref{thm:meanseg} of \S\ref{sec:maps_with_prop_p_main_thm} below where such maps are called maps with property $(\mathrm{P})$. It is shown that the closure of the mean action spectrum of maps with property $(\mathrm{P})$ contain a non-trivial interval. Finally we prove the last case where we assume:
\begin{enumerate}[resume]
\item $\phi: \D \rightarrow \D$ has zero topological entropy. Moreover, there is an iteration $\phi^{n}$ and an isotopy $\phi_{n,t}$ with time-one map $\phi^{n}$ satisfying $\vert \rho(\phi_{n,t})\vert \leq 1$ and either, 
\begin{enumerate}
\item $\sigma_{+}(\phi_{n,t})\neq \max\{-\rho(\phi_{n,t}), 0\}$ or $\sigma_{-}(\phi_{n,t}) \neq \min\{-\rho(\phi_{n,t}), 0\}$.
\end{enumerate}
\end{enumerate}
Let $\psi = \phi^{n}$ and $\psi_{t} = \phi_{n,t}$ for simplicity. Without loss of generality we can assume that $\sigma_{+}(\psi_{t})\neq \max\{-\rho(\psi_{t}),0\}$ as if $\sigma_{-}(\psi_{t}) \neq \min\{-\rho(\psi_{t}), 0\}$ then we would work with $\psi^{-1}$ and $\psi_{t}^{-1}$. Moreover, we can further reduce the problem to the case of $\sigma_{+}(\psi_{t}) > -\rho(\psi_{t})$ because the remaining case $-\rho(\psi_{t}) > \sigma_{+}(\psi_{t})$ can be reduced to the previous case by considering $\psi_{t}' = \psi_{t}^{-1}\circ R_{-2\pi t}$ which would imply,
$$\sigma_{+}(\psi_{t}') \geq \sigma_{+}(R_{-2\pi t}) - \sigma_{+}(\psi_{t}) > 1 + \rho(\psi_{t}) = -\rho(\psi_{t}').$$ 
We prove the following statement to finish the proof of Theorem \ref{thm:irratseg}. By the above discussion we are assuming $\sigma_{+}(\psi_{t}) > \max\{-\rho(\psi_{t}), 0\}$. In the following theorem, mean action of a periodic point $x$ of period $k$ for $G$ is meant to be $\mathcal{A}_{G^{(k)}}(x)/k = -a_{\psi_{G,t}^{k}}(x)/k$; see Remark \ref{rem:compare_two_action_defs}.
\begin{theorem}\label{thm:irratseg}
Let $\psi$ and $\psi_{t}$ be as above. Then, for every integer $k$ and prime number $p$ and any interval $[a,b] \subset (\max\{-\rho(\psi_{t}), 0\}, \sigma_{+}(\psi_{t})]$ satisfying,
$$p^{k} > \max\{\frac{1}{a}, \frac{1}{b-a}\},$$
there is an $p^{k}$-periodic point of $\psi$ with mean action inside the interval $[a,b]$.
\end{theorem}
\begin{proof} We first look at the case $\rho(\psi_{t})> 0.$ The numbers $\rho(\psi_{t})$ and $\sigma_{+}(\psi_{t})$ both remain unchanged if one switches to another isotopy $\psi_{t}'$ with time-one map $\psi$ that is homotopic to $\psi_{t}$ as long as the generating Hamiltonians vanish on the boundary $\partial \D$. Therefore, throughout the proof we select the isotopy $\psi_{t}$ given by Lemma \ref{lem:specialham}, whose generating Hamiltonian is denoted by $H$, and consider its extension $\psi_{\beta,t}$ given by Lemma \ref{lem:extcompact}; so $\psi_{\beta,1}$ is a compactly supported map inside a disk $\D_{1 + \beta}$ of radius $1 + \beta$; see also Figure \ref{fig:extended_hamilt}. by property $(3)$ of Theorem \ref{lem:extcompact} the action function of $\psi_{\beta,t}$ satisfies,
\begin{equation}\label{ineqrot}
0 \leq {a_{\psi_{\beta,t}}|}_{\D_{1 + \beta}\backslash \D} < {a_{\psi_{\beta,t}}|}_{\partial \D} = \rho(\psi_{t})+a_{\beta} 
\end{equation}
which for small enough $\beta$ implies that $\sigma_{+}(\psi_{\beta,t}) = \sigma_{+}(\psi_{t}) - a_{\beta}$, where $a_{\beta}>0$ is a small number converging to zero as $\beta \to 0$, see Lemma \ref{lem:extcompact} point (2) and Figure \ref{fig:extended_hamilt}, by spectrality and continuity of the extended spectral invariants; see Theorem \ref{T:spec main-1 Intro} and Remark \ref{rmk:easier_extension}. By Theorem \ref{thm:compmean}, for small enough $\beta$, we know that for the interval,
\begin{equation}\label{eq:shifted_int}
[a - a_{\beta}, b - a_{\beta}] \subset (\rho(\psi_{t}) - a_{\beta}, \sigma_{+}(\psi_{t}) - a_{\beta}] \subset (0, \sigma_{+}(\psi_{\beta,t})],
\end{equation}
and the number $k,p$ as in the statement there is an $p^{k}$-periodic point $x$ of $\psi_{\beta,1}$ with 
$$\mathcal{A}_{H_{\beta}^{(p^{k})}}(x)/p^{k} \in [a - a_{\beta}, b - a_{\beta}],$$
where $H_{\beta}$ is a compactly supported Hamiltonian generating $\psi_{\beta,t}$, see Remark \ref{rmk:easier_extension}; note that $a_{\beta} \to 0$ as $\beta \to 0$ and so if $p^{k} > 1/a$ then for small enough $\beta$ we have $p^{k} > 1/(a + a_{\beta})$. It follows from \eqref{ineqrot} and \eqref{eq:shifted_int} and the fact that $\mathcal{A}_{H_{\beta}} = -a_{\psi_{\beta,1}}$, that the detected periodic point $p$ must be inside $\D$, as periodic points outside of $\D$ have negative symplectic action, i.e., $\mathcal{A}_{H_{\beta}} \leq 0$ on $\D_{1 + 1/n}\backslash \D$. Hence, $x$ must be a periodic point of $\psi$. Since $a_{\psi_{\beta,t}^{r}} = a_{\psi_{t}^{r}} + ra_{\beta}$ on the disk $\D$, see item $(2)$ of Lemma \ref{lem:extcompact}, we conclude that,
$$\mathcal{A}_{H^{(p^{k})}}/p^{k} = \mathcal{A}_{H_{\beta}^{(p^{k})}}/p^{k} + a_{\beta}, \text{ which implies, } \mathcal{A}_{H^{(p^{k})}}(x)/p^{k} \in [a,b].$$ 
This finishes the proof in this case. 

Let us now assume that $\rho(\psi_{t}) < 0$ and $\sigma_{+}(\psi_t) >- \rho(\psi_t)$. We work with the isotopy $\psi_{t}$ given by Lemma \ref{lem:specialham} as before. We then consider $\psi_{t}^{-1}$ for which one has $\rho(\psi_{t}^{-1}) > 0$. We extend it to the compactly supported isotopy $\bar{\psi}_{\beta, t}$, generated by $H_{\beta}'$ as in Remark \ref{rmk:easier_extension} and Lemma \ref{lem:extcompact}. Let $\psi_{\beta,t} = (\bar{\psi}_{\beta,t})^{-1}$, generated by $H_{\beta} = \bar{H}_{\beta}'$; the latter isotopy will be a compactly supported extension of $\psi_{t}$; see Figure \ref{fig:extended_hamilt_2} below. 
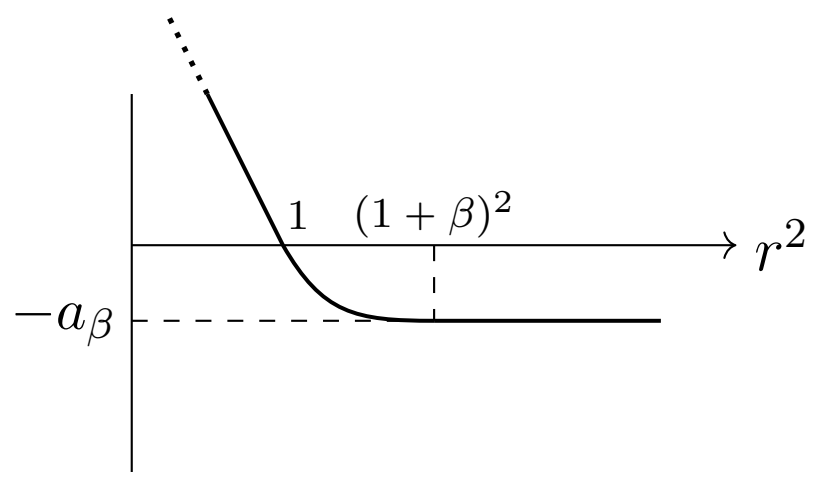
\begin{figure}[ht]
\centering
    \centering
    \begin{tikzpicture}[scale = 1]
    \draw (-1,-1.5) -- (-1,1);
    \draw[->] (-1,0) -- (3, 0);\node at (3.3, 0) {$r^{2}$};
    \draw[line width=0.75pt] (-0.5,1) -- (0,0);
    \draw[dotted][line width=1pt] (-0.75,1.5) -- (-0.5,1);
    \draw[line width=0.75pt] (1,-0.5) -- (2.5,-0.5);
    \draw[dashed] (-1,-0.5) -- (1,-0.5); 
    \draw[dashed] (1,0) -- (1,-0.5); 
    \node at (-1.45, -0.5) {\small{$-a_{\beta}$}};
    \node at (0.1, 0.2) {\tiny{$1$}};\node at (1, 0.2) {\tiny{$(1+\beta)^{2}$}};
    \draw[line width=0.75pt] (0,0)to[out=-60,in=180,looseness=1.2](1,-0.5);
    \end{tikzpicture}
\caption{The shifted extended Hamiltonian $H_{\beta} - a_{\beta}$.}
\label{fig:extended_hamilt_2}
\end{figure}
It is then not hard to check that the $\mathcal{A}_{H_{\beta}}$-action of one-periodic orbits inside $\D_{1+\beta}\backslash \D$ are positive and in the interval $[0, -\rho(\psi_{t}) + a_{\beta}]$. Take $[a,b] \subset (\max\{-\rho(\psi_{t}), 0\}, \sigma_{+}(\psi_{t})] = (-\rho(\psi_{t}), \sigma_{+}(\psi_{t})]$. For small enough $\beta$ we must have,
$$[a + a_{\beta}, b + a_{\beta}] \subset (-\rho(\psi_{t})+a_{\beta}, \sigma_{+}(\psi_{\beta,t})];$$
note that in this case $\sigma_{+}(\psi_{\beta,t}) = \sigma_{+}(\psi_{t}) + a_{\beta}$; where $a_{\beta} > 0$ is a small constant as before. Applying Theorem \ref{thm:compmean} we obtain an $p^{k}$-periodic point $x$ with, 
$$\mathcal{A}_{H_{\beta}^{(p^{k})}}(x)/p^{k} \in [a + a_{\beta}, b +a_{\beta}].$$ 
Since the action of points in $\D_{1+\beta}\backslash \D$ are in the interval $[0, -\rho(\psi_{t}) + a_{\beta}]$, the detected periodic point $x$ must lie inside $\D$ in which case we have,
$$\mathcal{A}_{H^{(p^{k})}}(x)/p^{k} = \mathcal{A}_{H_{\beta}^{(p^{k})}}(x)/p^{k}- a_{\beta},$$
which implies $\mathcal{A}_{H^{(p^{k})}}(x)/p^{k}\in [a,b]$ and this finishes the proof.
\end{proof}

\subsection{Maps with property \texorpdfstring{$(\mathrm{P})$}{}}\label{sec:maps_with_prop_p_main_thm}
Throughout this section we often work with restriction of symplectic map $\phi$ to an invariant open topological annulus $A$. The extension of these restrictions $\phi|_{A}$ to the prime end compactification $A_{c}$ of $A$ will always be isotopic to the identity by Remark \ref{rem:rest_are_isotop_to_id}. So, we will freely apply the results of \cite{franks1988generalizations, franks2013entropy} to $\phi|_{A}$ whenever we need them.

We now introduce a class of highly degenerate maps for which, at the level of spectral invariants, one cannot distinguish them from pseudo-rotations. This can occur only when the rotation on the boundary is irrational. In the following when we write $\rho(\phi)$ we shall suggestively mean the rotation number of $\phi$ restricted to the boundary $\partial \D$. We say that an area-preserving map $\phi:\mathbb{D}\to\mathbb{D}$ has property $(\mathrm{P})$ if the following properties hold:
\begin{enumerate}
    \item $\rho(\phi)\in \mathbb{R}\backslash \mathbb{Q}$,
    \item $\phi$ has zero topological entropy.
    \item For every $n$ and isotopy $\phi_{n,t}$ with time-one map $\phi^{n}$ and $\vert\rho(\phi_{n,t})\vert < 1$ we have, 
    $$\sigma_{+}(\phi_{n,t})= \max\left\{-\rho(\phi_{n,t}),0\right\}, \ \sigma_{-}(\phi_{n,t}) = \min\left\{-\rho(\phi_{n,t}),0\right\}.$$
\end{enumerate} 
Recall that $\sigma_{\pm}$ are the extended spectral invariants from \eqref{eq:spec_inv_non_comp_supp_maps}. All pseudo-rotations satisfy property $(\mathrm{P})$ above. In the following by a periodic topological disk $D$ for $\phi$ we mean that for some $k > 1$, $\phi^{j}(D) \cap D = \emptyset$ for $1 \leq j \leq k-1$ and $\phi^{k}(D) = D$.
\begin{lemma}\label{lem:noperdisc}
Maps with property $(\mathrm{P})$ do not have any periodic topological disks.
\end{lemma}
\begin{proof}
Assume there exists a periodic topological disc $D_{0}$ with period $k$, i.e.\ $\phi^{k}(D_{0}) = D_{0}$. 
Consider sequences $\{p_{n}\}$ and $\{q_{n}\}$ with $q_{n} \to \infty$ such that $k \nmid q_{n}$ for all $n$, and, 
$$|q_{n}\alpha - p_{n}| \leq 1/n;$$  
this is easy to see as one can choose such $p_{n}$ and $q_{n}$ with $\vert q_{n}\alpha - p_{n}\vert$ being arbitrarily small hence $n < 1/\vert q_{n}\alpha - p_{n}\vert.$ Choose an isotopy $\phi_{q_{n},t}$ with time-one map $\phi^{q_{n}}$ with boundary rotation number equal to $q_{n}\alpha - p_{n}$. By definition of the $\gamma$-norm, we then have
$$\gamma(\phi^{q_{n}}_{t}) = |q_{n}\alpha - p_{n}|.$$  
Since $k \nmid q_{n}$ it follows that $\phi^{q_{n}}$ displaces the disc $D_{0}$. Hence, by item $(4)$ of Theorem \ref{T:Vit 1},
$$1/n \geq \vert q_{n}\alpha - p_{n}\vert = \gamma(\phi^{q_{n}}_{t})\geq \mathrm{area}(D_{0})/2,$$  
which is a contradiction for large enough $n$. This completes the proof.
\end{proof}

Now we prove an important result about the existence of a ``foliation" near the boundary. In the following lemma an annulus $A'$ inside another annulus $A$ is called essential if it is non-contractible inside $A$.

\begin{lemma}\label{lem:link}
Suppose $\phi:\mathbb D\to\mathbb D$ satisfies property $(\mathrm{P})$ and is not a pseudo-rotation. Then, there exists an annulus $U\in\mathcal A(\phi)$ with one boundary component being $\partial\mathbb D$ and $U$ contains no fixed points of $\phi$. 
Moreover, for every integer $N$, there is a collection of at least $N$ disjoint, essential, regular and $\phi$-invariant open annuli inside $U$, so that each frontiers of these invariant annuli contains at least one periodic point of $\phi$.
\end{lemma}
\begin{proof}
Fix a choice of an isotopy $\phi_{t}$ with time-one map $\phi$ throughout the proof.  

Step $1$. By the Brouwer plane translation theorem \cite{MR1511684}, see also the version on orientation-preserving homeomorphisms of the plane with a non-wandering point by Fathi \cite{MR925994}, there is a fixed point $x_{0}$ of $\phi$ in the interior of the disk 

Step $2$. Since property $(\mathrm{P})$, and being a pseudo-rotation, are invariant under conjugation by compactly supported area preserving maps, we may assume that $x_{0}$ is the origin.

Step $3$. Collapse the boundary of the disk to a point $p$ to obtain a sphere $S$ and let $\bar{\phi}:S\rightarrow S$ be the area-preserving homeomorphism induced by $\phi$ which fixes the point $p$. Applying Theorem \ref{thm:franks_handel_theorem} to $\bar{\phi}$, there is an open annulus $U''$ one of whose frontiers is the single point $p$. Let 
\begin{equation}\label{eq:rot_numb_ann}
U' := U'' \backslash \{p\} \subset \D\backslash \partial \D
\end{equation}
that covers a neighborhood of the boundary and contains no fixed point of $\phi$; notice that the map $\phi$ has no fixed points near the boundary and hence the collapsed boundary in the sphere is an isolated fixed point. Note that $U$ is essential in $\D^{\times} = \D \backslash\{0\}$.

Step $4$. We show that up to passing to an iteration of $\phi$ one can find another annulus $U\Subset U'$ where one can assume that rotation numbers on $U$ and \emph{on its frontiers} are all well-defined. It is shown in \cite[Theorem 1.4(1)]{franks2013entropy} that the rotation number of any point inside $U'$, potentially not on its frontiers, is well-defined and since $\phi$ has zero topological entropy the rotation number is also continuous; the continuity is not used in the rest of the proof. As $\phi$ is not a pseudo-rotation, it follows from \cite[Theorem 1.4(2)]{franks2013entropy} that $\rho_{\widetilde{\phi}_{c}}$ is not constant on $U'$ where $\widetilde{\phi}_{c}$ is a lift $(\phi|_{U'})_{c}$ to the universal cover $\widetilde{U}'_{c}$ where $(\phi|_{U'})_{c}:U'_{c}\rightarrow U'_{c}$ is the prime extension of $\phi|_{U'}$ to the prime end compactification of $U'$. Hence, the rotation interval $\widetilde{\phi}_{c}(U'_{c})$, see \cite[Theorem 2.3]{franks2013entropy} and \cite{Handel-90-rot-set-is-closed-int}, is non-trivial and we can apply \cite[Lemma 15.9]{franks2013entropy} to find a large prime number $q$ and some $V\in \mathcal{A}(q)$ so that $\bar{V} \subset U'$. Hence, up to passing to the iteration $\phi^{q}$, we find an invariant annulus 
\begin{equation}\label{eq:the_ann_U}
U \text{ whose frontiers are the interior frontier of V and } \partial \D.
\end{equation}
So, we can assume that the rotation number is well-defined on $\bar{U}$ defined by
\begin{equation}\label{eq:the_univ_rot_numb}
\rho': \bar{U}\rightarrow \R; z\mapsto \rho_{\widetilde{\phi}_{c}}(z).
\end{equation}
The rotation number on the interior frontier of $U$ is an integer as it contains a fixed point of $\phi$; see Theorem \ref{thm:franks_handel_theorem} above.

Step 5. In the previous step we found an \emph{essential} annulus $\bar{V}\subset U'$, by only appealing to \cite[Lemma 15.9]{franks2013entropy} without using Property $(\mathrm{P})$, and used it to find $U$ on closure of which the rotation number is well-defined everywhere. Now we use property $(\mathrm{P})$ to prove that all annuli $V\in \mathcal{A}(q)$ that intersect $U$ are \emph{essential} in $U$, for all $q\geq 2$. Fix an integer $q \geq 2$. Consider $V \in \mathcal{A}(q)$ that intersects $U$; of course such an annulus exists. By item $(3)$ of Theorem \ref{thm:regopen} we must have $V \subset U$. By item $(2)$ of Theorem \ref{thm:regopen}, for every annulus $V\in \mathcal{A}(q)$ for $q\geq 2$, it must be that either $\phi(V)=V$ or $\phi$ permutes $V$ within $\mathcal{A}(q)$, that is there exists some $k>1$ so that,
$$\phi^{j}(V) \cap V = \emptyset, 1\leq j\leq k-1, \phi^{k}(V)=V.$$ 
The latter cannot happen by Lemma \ref{lem:noperdisc}. Hence, all annuli of $\mathcal{A}(q)$ are fixed by $\phi$. It then follows that any annulus $V \in \mathcal{A}(q)$ that lies inside $U$ must be essential in $\mathbb{D}^{\times}$. Indeed, if $V$ were not essential and $\phi(V)=V$, then the outer frontier of $V$ would bound an open disk $D_{V}$ in $U$ such that $\phi(D_{V})=D_{V}$. By the Brouwer plane translation theorem \cite{MR1511684} and the fact that $\phi$ is area-preserving, see also \cite{MR925994}, $\phi$ would then have a fixed point in the interior of $D_{V}\subset U$, which is impossible since $\phi$ has no fixed points in the interior of $U$. Hence, all the annuli $V \in \mathcal{A}(q)$ that intersect $U$ are essential in $U$ and are $\phi$-invariant for all $q\geq 2$, and each frontier of these annuli contain at least one periodic point of $\phi$. By item $(1)$ of Theorem \ref{thm:regopen} all of these annuli are open and regular.

Step $6$. We now prove that for every $N\geq 1$ there are at least $N$ disjoint such invariant annuli as in Step $5$. This basically follows from \cite[Theorem 2.3(1)]{franks2013entropy} where it is shown that any rational number in the rotation interval $\rho'(U)$ is realized by a periodic point; the function $\rho'$ is defined in \eqref{eq:the_univ_rot_numb} and here we use \cite[Lemma 15.7(1)]{franks2013entropy} which implies that the rotation number defined via $(\phi|_{U})_{c}$ is equal to $\rho'$ since $U$ is essential in $U'$. Namely, consider distinct rational numbers $p_{i}/q_{i}\in \rho'(U), i = 1,\dots, N,$ with $q_{i}$'s being distinct prime numbers and the corresponding periodic points $x_{i}, i = 1,\dots, N$, then, for a large enough prime number $q > \max q_{i}$, by Lemma \ref{rem:rotbound} the points $x_{i}$ must lie in distinct elements $V_{i}\in \mathcal{A}(q)$. This shows the existence of at least $N$ disjoint $\phi$-invariant $V_{i} \subset U$ belonging to $\mathcal{A}(q)$ for some $q$ and they have periodic points of $\phi$ on their frontiers. 
\end{proof}

\begin{lemma} \label{rem:rotbound}
Let $U$ be as constructed in the proof of Lemma \ref{lem:link}. For any $V \in \mathcal{A}(q)$ with $V \subset U$, and every $x, y \in \bar{V}$ we have the following inequality: 
$$\vert\rho'(x)-\rho'(y)\vert \leq 1/q,$$
where $\rho': \bar{U}\rightarrow \R$ is defined in \eqref{eq:the_univ_rot_numb}.
\end{lemma}
\begin{proof}
For such $V$, there are no interior $q$-periodic points of $\phi$, and since $V$ is essential in $U$ and hence in $U'$, see \eqref{eq:rot_numb_ann} and \eqref{eq:the_ann_U}, the rotation number $\rho_{V}:V_{c}\rightarrow \R$ defined via $(\phi|_{V})_{c}$ coincides with $\rho'$ on $V$ modulo $\Z$; see \cite[Lemma 15.7(1)]{franks2013entropy}. Suppose 
$\vert\rho'(x)-\rho'(y)| >1/q,$ for some points $x,y\in \bar{V}$. Since $\rho'$ is continuous we can assume that $x,y\in V$. Then, it follows that
$$\vert\rho_{V,q}(x)-\rho_{V,q}(y)| > 1,$$ 
where $\rho_{V,q}; V_{c}\rightarrow \R$ is defined via $(\phi^{q}|_{V})_{c}$ and it equals $\rho_{V,q} = q\rho_{V}$. This together with Franks' theorem \cite{franks1988generalizations} contradicts the fixed point free-ness of $\phi^{q}$ in $V$; the rotation interval of every lift of $\phi^{q}$ will have an integer in its interior, so one can find a lift of $\phi^{q}$ to the universal cover of $V_{c}$ with the origin in the interior of its rotation interval, this results in both positively and negatively recurrent disks for the lift and one can apply \cite{franks1988generalizations} to find a fixed point in the annulus $V$. 
\end{proof}

\begin{lemma} \label{lem:segrot}
Let $V$ be an open regular, $\phi$-invariant, and essential annulus with $\bar{V}\subset U$. Fix a continuous curve $\gamma:[0,1]\to \bar{V}$ with $\gamma(0)$ and $\gamma(1)$ lying on different frontiers of $V$. Then, we  have $\rho'(\bar{V})=\rho'(\gamma([0,1]))$ where $\rho'$ is defined in \eqref{eq:the_univ_rot_numb}.
\end{lemma}

\begin{proof}
The inclusion $\rho'(\gamma([0,1])) \subseteq \rho'(\bar{V})$ is obvious.
We prove the opposite inclusion $\rho'(\bar{V}) \subseteq \rho'(\gamma([0,1])).$ Assume that there exists $a \in \rho'(\bar{V}) \backslash \rho'(\gamma([0,1])).$ Since both $\rho'(\gamma([0,1]))$ and $\rho'(\bar{V})$ are compact connected sets, there exists $\epsilon > 0$ such that 
$$[a-\epsilon, a+\epsilon] \cap \rho'(\gamma([0,1])) = \emptyset.$$
Choose a rational number $p'/q' \in (a - \epsilon, a + \epsilon)$ for some prime number $q'$. By \cite[Theorem 2.3]{franks2013entropy} there is a periodic point $x\in \bar{V}$ with period $q'$ and rotation number $p'/q'.$ In fact $x \in V$ since the rotation numbers of the frontiers are included in $\rho'(\gamma([0,1]))$. For a large enough prime number $q'' > q'$ there will be an annulus $V''\in \mathcal{A}(q'')$ such that $x \in V''$. As shown in the proof of Lemma \ref{lem:link} the annulus $V''$ must be essential in $U$; this also follows from \cite[Lemma 15.7(2)]{franks2013entropy} without the property $(\mathrm{P})$, because otherwise almost all the rotation numbers, including rotation numbers of periodic points inside $V''$, would be the same and equal to $p''/q''$ for some integer $0 < p'' < q''$ which is not the case for the periodic point $x$. This however is a contradiction because by Lemma \ref{rem:rotbound} for large enough $q''$ we have $\rho'(V'') \subset [a - \epsilon, a + \epsilon],$ and as $V''$ is essential we have $V'' \cap \gamma([0,1]) \neq \emptyset$ which implies that $\rho'(\gamma([0,1]))$ must intersect $[a - \epsilon, a + \epsilon]$.
\end{proof}

In the following claim $U$ is constructed in Lemma \ref{lem:link}; see \eqref{eq:the_ann_U}. Recall that by construction the closure of $U$ is contained in the interior of $U'$, see \eqref{eq:rot_numb_ann}, where the rotation number is defined. Hence, the rotation number is well-defined on the closure of $U$ as well via \eqref{eq:the_univ_rot_numb}.

\begin{claim}\label{clm:small_rot_small_mean_act}
Let $U$ be as constructed in Lemma \ref{lem:link} and $V \in \mathcal{A}(\phi,q)$, for some $q$, be an essential annulus inside $U$ such that the. Let $\phi_{t}$ be a symplectic isotopy with time-one map $\phi$. Let $\rho':\bar{U}\rightarrow \R$ be as in \eqref{eq:the_univ_rot_numb}. Assume that for some $\epsilon >0$ and every two points $x,y\in \bar{V}$ one has $$\vert\rho'(x)-\rho'(y)\vert\leq \epsilon.$$ 
Then, for any two periodic points $x_{0},x_{1}$ lying on different frontiers $F_{0},F_{1}$ of $V$ we have 
$$\vert a^{\infty}_{\phi_{t}}(x_{1})-a^{\infty}_{\phi_{t}}(x_{0})\vert\leq \epsilon.$$ 
\end{claim} 
\begin{proof}
Let $q > 1$ be an integer so that $x_{0}, x_{1} \in \Fix(\phi^{q})$ which in particular implies $a^{\infty}_{\phi_{t}^{}}(x_{i}) = a_{\phi_{t}^{q}}(x_{i})/q.$ Let $D_{0}, D_{1}$ be the bounding disks of the interior frontier $F_{0}$ and exterior frontier $F_{1}$ of $V$ respectively. 

Claim. For every point $y\in \mathrm{int}(D_{i})$ one has $W^{\infty}_{\phi_{t}}(x_{i}, y) = W^{\infty}_{\phi_{t}}(x_{i}, 0)$. Similarly, for every point $y\in \D \backslash D_{i}$ one has $W^{\infty}_{\phi_{t}}(x_{i},y) = W^{\infty}_{\phi_{t}}(0, y)$.

Proof of Claim. The first statement is clear. The second statement is also easily shown for $y \in \D\backslash \bar{D}_{i}$. The non-trivial part is to show that for $y\in F_{i} = \bar{D}_{i}\backslash D_{i}$ one has 
$$W^{\infty}_{\phi_{t}}(x_{i},y) = W^{\infty}_{\phi_{t}}(0, y).$$
To see this note that $W^{\infty}_{\phi_{t}}(x_{i},z)$ and $W^{\infty}_{\phi_{t}}(0,z)$ are equal to some rotation numbers and consequently they are continuous functions in $z$. If $y\in F_{i}$ is not a $q$-periodic point of $\phi$ then by regularity of Franks-Handel invariant annuli, i.e., $\mathrm{int}(\bar{V}) = V\subset \D\backslash \Fix(\phi^{q})$, one can find a sequence of points $\{y_{k}\}\subset \D\backslash \bar{D}_{i}$ converging to $y$. Hence,
$$W^{\infty}_{\phi_{t}}(x_{i},y)=\lim_{k\to\infty}W^{\infty}_{\phi_{t}}(x_{i},y_{k}) = \lim_{k\to \infty}W^{\infty}_{\phi_{t}}(0,y_{k}) = W^{\infty}_{\phi_{t}}(0,y).$$
Suppose now that $y$ is a $q$-periodic point. One can show that the function $W^{\infty}_{\phi_{t}}(w,z)$ is constant in $z\in F_{i}$ for any given $q$-periodic point $w\in F_{i}$: because for every $z\in F_{i}$ and any point $z^{*}\in \Fix(\phi^{q})$ in the $\omega$-limit set of $z$ it follows from Lemma \ref{lem:wind_omega_limit} that
$$W^{\infty}_{\phi_{t}}(w,z) = W^{\infty}_{\phi_{t}}(w,z^{*})\in \Z/q.$$
So, the function $W^{\infty}_{\phi_{t}}(w,z)$ is continuous in $z\in F_{i}$ and is valued in a discrete set and so it must be constant as $F_{i}$ is connected. Let $z^{*}\in F_{i}$ be any point that can be approximated by points $\{z_{k}\}\subset D_{i}$. Then we have,
$$\begin{aligned}
W^{\infty}_{\phi_{t}}(x_{i},y)&= W^{\infty}_{\phi_{t}}(y,x_{i})=W^{\infty}_{\phi_{t}}(y,z^{*})\\
&=\lim_{k\to \infty}W^{\infty}_{\phi_{t}}(y,z_{k}) =\lim_{k\to \infty}W^{\infty}_{\phi_{t}}(z_{k},y)\\
&=\lim_{k\to \infty}W^{\infty}_{\phi_{t}}(0,y)=W^{\infty}_{\phi_{t}}(0,y);
\end{aligned}$$
this finishes the proof of the claim. \qed

\bigskip

By Theorem \ref{thm:actwind}, we now obtain the following identity:
\begin{equation}\label{eq:mean_act_decomp}
a^{\infty}_{\phi_{t}}(x_{i})= \int_{\mathbb{D}}W^{\infty}_{\phi_{t}}(x_{i}, y)\,\omega= \int_{D_{i}}W^{\infty}_{\phi_{t}}(x_{i}, 0)\omega + \int_{\mathbb{D}\backslash D_{i}}W^{\infty}_{\phi_{t}}(0, y)\omega,
\end{equation}
for $i = 0,1$. Denote the area of $D_{i}$ by $a_{i}$. So,
$$\begin{aligned}
\big\vert a^{\infty}_{\phi_{t}}(x_{1}) - a^{\infty}_{\phi_{t}}(x_{0})\big\vert &= \bigg\vert a_{1}W^{\infty}_{\phi_{t}}(x_{1}, 0) - a_{0}W^{\infty}_{\phi_{t}}(x_{0}, 0) - \int_{D_{1}\backslash D_{0}}W^{\infty}_{\phi_{t}}(0, y)\omega\bigg\vert\\
&=\bigg\vert a_{0}\big(W^{\infty}_{\phi_{t}}(x_{1}, 0) - W^{\infty}_{\phi_{t}}(x_{0}, 0)\big) - \int_{D_{1}\backslash D_{0}}\big(W^{\infty}_{\phi_{t}}(x_{1}, 0) - W^{\infty}_{\phi_{t}}(0, y)\big)\omega\bigg\vert\\
&\leq a_{0}\big\vert\rho'(x_{1}) - \rho'(x_{0})\big\vert + \int_{D_{1}\backslash D_{0}}\big\vert\rho'(x_{1}) - \rho'(y)\big\vert\omega\\
&\leq a_{0}\epsilon + (a_{1} - a_{0})\epsilon \leq \epsilon,
\end{aligned}$$
which finishes the proof. On the third line we used $W^{\infty}_{\phi_{t}}(x,0) = \rho'(x)$ for every $x\in \bar{V}$.
\end{proof}
\begin{lemma}\label{lem:wind_omega_limit}
Let $V\subset U$ be as in Claim \ref{clm:small_rot_small_mean_act}. Let $x_{i}\in F_{i}, i=0,1,$ be $q$-periodic points of $\phi$ on the frontiers of $V$. Let $z\in F_{i}$ be any point and $z^{*}\in F_{i}\cap \Fix(\phi^{q})$ be a point in the $\omega$-limit set of $z$. Then,
$$W^{\infty}_{\phi_{t}}(x_{i}, z) = W^{\infty}_{\phi_{t}}(x_{i}, z^{*}).$$
\end{lemma}
\begin{proof}
Let $\psi = \phi^{q}$ and $\psi_{t} = \phi_{t}^{q}$, and observe that $\psi(x_{i}) = x_{i}$. Then, we have,
$$W^{\infty}_{\psi_{t}}(x_{i}, \psi(z))= W^{\infty}_{\psi_{t}}(x_{i}, z).$$
Therefore, if $\phi^{n_{k}}(z)\to z^{*}$ for some infinite sequence $\{n_{k}\}_{k\geq 1}$, it follows that,
$$W^{\infty}_{\phi_{t}}(x_{i},z)=\frac{1}{q}W^{\infty}_{\psi_{t}}(x_{i},z)=\frac1q\lim_{k\to\infty} W^{\infty}_{\psi_{t}}(x_{i}, \psi^{n_{k}}(z)) = \frac1q W^{\infty}_{\psi_{t}}(x_{i}, z^{*})=W^{\infty}_{\phi_{t}}(x_{i}, z^{*}),$$
which finishes the proof.
\end{proof}
Now we prove Theorem \ref{mainthm:meanacion} for maps with property (P).
\begin{theorem} \label{thm:meanseg}
Suppose $\phi: \D \rightarrow \D$ has property $(\mathrm{P})$ and is not a pseudo-rotation. Then, the closure of the mean action spectrum of $\phi$ contains a non-trivial interval.
\end{theorem}
\begin{proof}
Fix a symplectic isotopy $\phi_{t}$ with time-one map $\phi$. Let $U$ be the invariant annulus from Lemma \ref{lem:link} defined in \eqref{eq:the_ann_U}, where we can assume that the rotation number $\rho':\bar{U}\rightarrow \R$ is defined on $\bar{U}$ and it is continuous, consequently, the rotation interval $\rho'(\bar{U})$ is non-trivial, see also \cite{Handel-90-rot-set-is-closed-int}; note that the rotation number on the boundary $\partial \D$ is an irrational number and it is a rational number on the interior frontier of $U$ so the rotation interval must be non-trivial as it contains at least two distinct numbers. Assume that 
\begin{equation}\label{eq:the_per_point}
\exists x\in \bar{U} \text{ such that } a^{\infty}_{\phi_{t}}(x) \neq \rho(\phi_{t}),
\end{equation}
where the latter is the rotation number of the isotopy $\phi_{t}$ on the boundary $\partial \D$. Up to passing to an iteration of $\phi$ and replacing U with another essential annulus in $U$ we might assume that $x$ is on the interior frontier of $U$. Define the following (mean action) interval:
\begin{enumerate}
\item $\mathrm{MAI}:=[a^{\infty}_{\phi_{t}}(x), \rho(\phi_{t})]$ if $\rho(\phi_{t}) > a^{\infty}_{\phi_{t}}(x)$; work with $\mathrm{MAI}=[\rho(\phi_{t}), a^{\infty}_{\phi_{t}}(x)]$ otherwise.
\end{enumerate}
Let $q$ be an arbitrary large number. Since the rotation number $\rho':\bar{U}\rightarrow \R$ is continuous on $\bar{U}$ there exists $\epsilon > 0$ so that for every two points $z,w\in \bar{U}$ with $d(z,w) < \epsilon$ one has
\begin{equation}\label{eq:small_dist_small_rot}
|\rho'(z) - \rho'(w)| \leq 1/q.
\end{equation}
Let $V_{1}, \dots, V_{N}\in \mathcal{A}(q)$ be the essential open regular annuli that are inside $U$; these are all shown to be essential inside $U$, see Lemma \ref{lem:link}, see also \cite[Theorem 1.2]{franks2013entropy} or Theorem \ref{thm:franks_handel_theorem} above. 

Pick a smooth simple curve $\gamma: [0,1] \rightarrow \overline{U}$ with finite length so that $\gamma(0)\in \mathrm{fr}_{0}(U)$ and $\gamma(1)\in \mathrm{fr}_{1}(U) = \partial \D$ and $\gamma((0,1)) \subset U$; such a curve exists basically because there exists a bi-holomorphic $F:\A\rightarrow U$ from the standard annulus $\A=\{z\in \C: a<|z|<b\}$ for some numbers $0<a<b<\infty$, and the energy of the bi-holomorphism $F$ is finite. So one can take the simple curve $\Gamma:(0,1)\rightarrow \A$ which is a re-parametrization of  $\{(r,\theta): \theta =\theta_{0}, a<r<b\}$ that has finite length connecting the two boundaries of $\A$ and let $\gamma:= F\circ \Gamma$, then $\gamma$ will have finite length for a generic $\theta_{0}$; see for instance \cite{MR2193309}. Therefore, 
$$\lim_{r\to 0}F\circ\Gamma(r) \text{ and } \lim_{r\to 1}F\circ\Gamma(r)$$
exist and lie on $fr_{0}(U)$ and $fr_{1}(U)=\partial \D$ respectively. Let $x_{0} = \gamma(0)\in fr_{0}(U)$ and $x_{k} = \gamma(1)\in \partial\D$. Let $x_{i} = F\circ \Gamma(\frac{i}{k}) \in \mathrm{img}(\gamma), i = 1, \dots, k-1$ be an ordered set of points that are ``equally'' distanced within $\mathrm{img}(\gamma)$. It follows from \eqref{eq:mean_act_decomp} that
$$a^{\infty}_{\phi_{t}}(x_{0}) = a^{\infty}_{\phi_{t}}(x).$$
Choose $k$ large enough so that $\mathrm{diam}([x_{i}, x_{i+1}])< \epsilon$ for all $i$ where $[x_{i}, x_{i+1}]$ is the segment of $\gamma$ connecting $x_{i}$ to $x_{i+1}$. Let $y_{i}\in \mathrm{img}(\gamma)$ with $d(x_{i}, y_{i})$ small enough so that $\mathrm{diam}([y_{i}, y_{i+1}]) < \epsilon$ for all $i$ and moreover $y_{i}\in V_{n_{i}}$ for some $n_{i}\in \{1,\dots, N\}$; these points exist because the annuli from $\mathcal{A}(q)$ form the set of weakly free disk recurrent points which are open and dense in $\D$ hence the union of the $V_i$ is open and dense in $U$. Pick $q$-periodic points $p_{i,0}, p_{i,1}$ on the frontiers $\mathrm{fr}_{0}(V_{n_{i}}), \mathrm{fr}_{1}(V_{n_{i}})$ respectively. Let $U_{i}$ be the invariant annulus with frontiers $\mathrm{fr}_{0}(U_{i}) := \mathrm{fr}_{1}(V_{n_{i}})$ and $\mathrm{fr}_{1}(U_{i}) := \mathrm{fr}_{0}(V_{n_{i+1}})$; see Figure \ref{fig:periodic_points_with_near_asymp_act}. It follows from \eqref{eq:small_dist_small_rot} and Lemma \ref{lem:segrot} applied to $U_{i}$ that for every two points $x, y\in U_{i}$ we have, 
$$|\rho'(x) - \rho'(y)| \leq \max\{|\rho'(z) - \rho'(w)|: z,w\in [y_{i}, y_{i+1}]\} \leq 1/q.$$
We now apply Lemma \ref{clm:small_rot_small_mean_act} to prove that,
$$|a^{\infty}_{\phi_{t}}(p_{i,1}) - a^{\infty}_{\phi_{t}}(p_{i+1,0})| \leq 1/q.$$
It also follows from Lemma \ref{rem:rotbound} and Lemma \ref{clm:small_rot_small_mean_act} that,
$$|a^{\infty}_{\phi_{t}}(p_{i,0}) - a^{\infty}_{\phi_{t}}(p_{i,1})| \leq 1/q,$$
for all $i = 1,\dots, k$, and one concludes the proof by choosing arbitrary large $q$'s, i.e., for every large number $q$ we found a collection of periodic points $\{p_{i,0}, p_{i,1}\}_{i=0}^{k}$ so that the biggest gap in the set 
$$\{a^{\infty}_{\phi_{t}}(p_{i,0}), a^{\infty}_{\phi_{t}}(p_{i,1})\}_{i=0}^{k} \cup\{a^{\infty}_{\phi_{t}}(x_{0}), a^{\infty}_{\phi_{t}}(x_{k})\}$$
is at most $1/q$ and this set contains the two ends of the interval $\mathrm{MAI}$, namely, $a^{\infty}_{\phi_{t}}(x_{0}) = a^{\infty}_{\phi_{t}}(x)$ and $a^{\infty}_{\phi_{t}}(x_{k}) = \rho(\phi_{t})$, so the set of mean actions of periodic points is dense in $\mathrm{MAI}$.
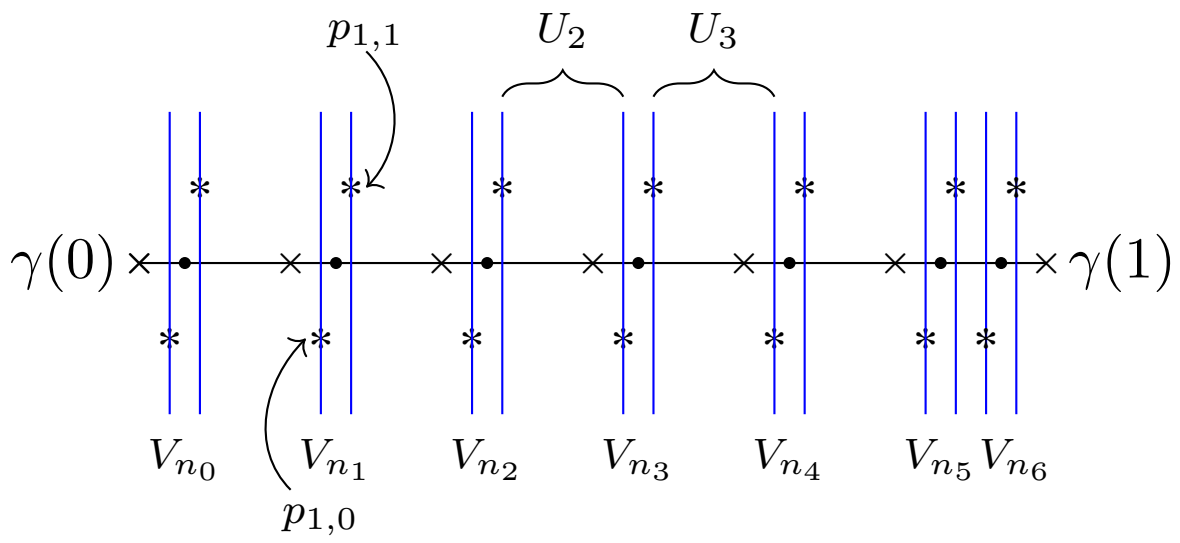
\begin{figure}[h]
  \centering
  \begin{tikzpicture}[scale = 1]
\draw (0,0) -- (6,0);
\node at (0, 0) {\tiny{$\times$}};
\node at (6, 0) {\tiny{$\times$}};
\node at (6.5, 0) {\small{$\gamma(1)$}};
\node at (-0.5, 0) {\small{$\gamma(0)$}};
\foreach \x in {0, 1, 2, ..., 5}{
    \node at (\x, 0) {\tiny{$\times$}};
    \node[circle, fill, inner sep=0.8pt] at (\x + 0.3, 0) {};
    \node at (\x + 0.4, 0.5) {\tiny{$*$}};
    \node at (\x + 0.2, -0.5) {\tiny{$*$}};
    \draw[color = blue] (\x + 0.2, -1) -- (\x + 0.2, 1);
    \draw[color = blue] (\x + 0.4, -1) -- (\x + 0.4, 1);
    \node at (\x + 0.3, -1.3) {\tiny{$V_{n_{\x}}$}};
}
\node[circle, fill, inner sep=0.8pt] at (6 - 0.3, 0) {};
\draw[color = blue] (6 - 0.4, -1) -- (6 - 0.4, 1);
\draw[color = blue] (6 - 0.2, -1) -- (6 - 0.2, 1);
\node at (6 - 0.2, -1.3) {\tiny{$V_{n_{6}}$}};
\node at (6 - 0.2, 0.5) {\tiny{$*$}};
\node at (6 - 0.4, -0.5) {\tiny{$*$}};
\draw [decorate, decoration={brace, amplitude=5pt}] (2.4,1.1) -- (3.2,1.1) 
        node [midway, above=5pt] {\tiny{$U_{2}$}};
\draw [decorate, decoration={brace, amplitude=5pt}] (3.4,1.1) -- (4.2,1.1) 
        node [midway, above=5pt] {\tiny{$U_{3}$}};

\draw[->] (1.5, 1.4) to[out = -45, in = 45, looseness = 1] (1.5, 0.5);
\draw[->] (1, -1.5) to[out = 135, in = -135, looseness = 1] (1.1, -0.5);
\node at (1 + 0.5, 1.5) {\tiny{$p_{1,1}$}};
\node at (1.2, -1.7) {\tiny{$p_{1,0}$}};
\end{tikzpicture}
\caption{The $x_{i}$'s are shown by $\times$, the $y_{i}$'s are shown by dots, and the periodic points $p_{i,0}, p_{i,1}$ are shown by $*$.}
\phantomsection\label{fig:periodic_points_with_near_asymp_act}
\end{figure}
It remains to prove the following claim:

Claim. There exists a periodic point $x$ satisfying \eqref{eq:the_per_point}. 

Proof of Claim. Assume that $V \in \mathcal{A}(q)$ and $V \subset U$ for some number $q>1$. We prove that any two periodic points lying on different frontiers of $V$ have different actions. This is sufficient, since it implies that at least one of them has action different from the boundary rotation number $\rho(\phi_{t})$. Let $x_{1}$ and $x_{2}$ be $q$-periodic orbits of $\phi$ lying on the frontiers $\mathrm{fr}_{0}(V)$ and $\mathrm{fr}_{1}(V)$ of $V$, respectively. Here, $\mathrm{fr}_{1}(V)$ is closer to the boundary of $\mathbb{D}$ than $\mathrm{fr}_{0}(V)$. We look at rotation numbers of $x_{1}$ and $x_{2}$. We have two cases: either $\rho'(x_{1})=\rho'(x_{2})\in \frac{1}{q} \mathbb{Z}$ or $\rho'(x_{1})\neq\rho'(x_{2})$ with $\rho'(x_{1}),\rho'(x_{2})\in \frac{1}{q} \mathbb{Z}$. In the former case, by the Franks fixed point theorem $\rho'$ cannot be the constant $\rho'(x_{1})$ on $V$ because $\phi^{q}$ has no fixed points in $V$. So, for every $x\in V$ we either have $\rho'(x)\geq \rho'(x_{1})$ and $\int_{V}\rho'(x)\omega>\int_{V}\rho'(x_{1})\omega$ or $\rho'(x)\leq \rho'(x_{1})$ and $\int_{V}\rho'(x)\omega < \int_{V}\rho'(x_{1})\omega$. Without loss of generality we assume the former case. Now, we compute the mean action of $x_{1}$ and $x_{2}$ using the formula \eqref{eq:mean_act_decomp} above:
$$\begin{aligned}
a^{\infty}_{\phi_{t}}(x_{1})-a^{\infty}_{\phi_{t}}(x_{2})=\int_{V}\rho'(y)\omega - \int_{V}\rho'(x_{1})\omega >0.
\end{aligned}$$ Hence the actions are different. We now look at the case $\rho'(x_{1})\neq\rho'(x_{2})$. Again, without loss of generality we assume $\rho'(x_{1})>\rho'(x_{2})$. By the Franks fixed point theorem, we have $\rho'(x_{1})\geq \rho'(x)\geq \rho'(x_{2})$ for every $x\in V$ and
$$\int_{V}\rho'(x_{1})\omega>\int_{V}\rho'(x)\omega> \int_{V}\rho'(x_{2})\omega.$$
We compute the mean action of $x_{1}$ and $x_{2}$ as before using \eqref{eq:mean_act_decomp},
$$\begin{aligned}
a^{\infty}_{\phi_{t}}(x_{1})-a^{\infty}_{\phi_{t}}(x_{2}) &= \int_{D_{1}}\rho'(x_{1})\omega -\int_{D_{2}}\rho'(x_{2})\omega+\int_{V}\rho'(x)\omega\\
&> \int_{D_{1}}\rho'(x_{2})\omega - \int_{D_{2}}\rho'(x_{2})\omega+\int_{V}\rho'(x)\omega\\
&= \int_{V}\rho'(x)\omega - \int_{V}\rho'(x_{2})\omega > 0.
\end{aligned}$$
Hence the actions are again different and this completes the proof of the theorem.
\end{proof}

\appendix
\section{Boundary depth is well-defined}\label{app:bound_depth_well_defined}
\begin{proof}[Proof of Proposition \ref{clm:ind_bound_depth}] 
Let $\delta' > \delta$ be a constant \emph{close} to $\delta$.  Let $\psi_{\epsilon}$ be a non-degenerate perturbation of $\phi \circ \phi^{1}_{\epsilon, \delta}$ generated by 
\begin{equation}\label{eq:gen_ham_eps}
H \# H_{\epsilon/2,\delta} \# H'_{\epsilon/2,\delta'} \# K,
\end{equation}
where $K$ is compactly supported inside $\{|z|^{2} < 1 - 2\delta\}$ and $H'_{\epsilon/2,\delta'}$ is a slight modification of $H_{\epsilon/2,\delta'}$; this is described in Figure \ref{fig:pert}. 
The purpose of this particular perturbation is to make the function $H \# H_{\epsilon,\delta}$, which is the Hamiltonian generating $\phi \circ \phi^{1}_{\epsilon, \delta}$, radial with non-zero slope near $|z|^{2} = 1-2\delta$ while keeping the slope $\epsilon$ at infinity, which enables us to take $K$ being supported in $\{|z|^{2} < 1-2\delta\}$ and that is why we break $H \# H_{\epsilon,\delta}$ into $H \# H_{\epsilon/2,\delta} \# H'_{\epsilon/2,\delta'}$. 

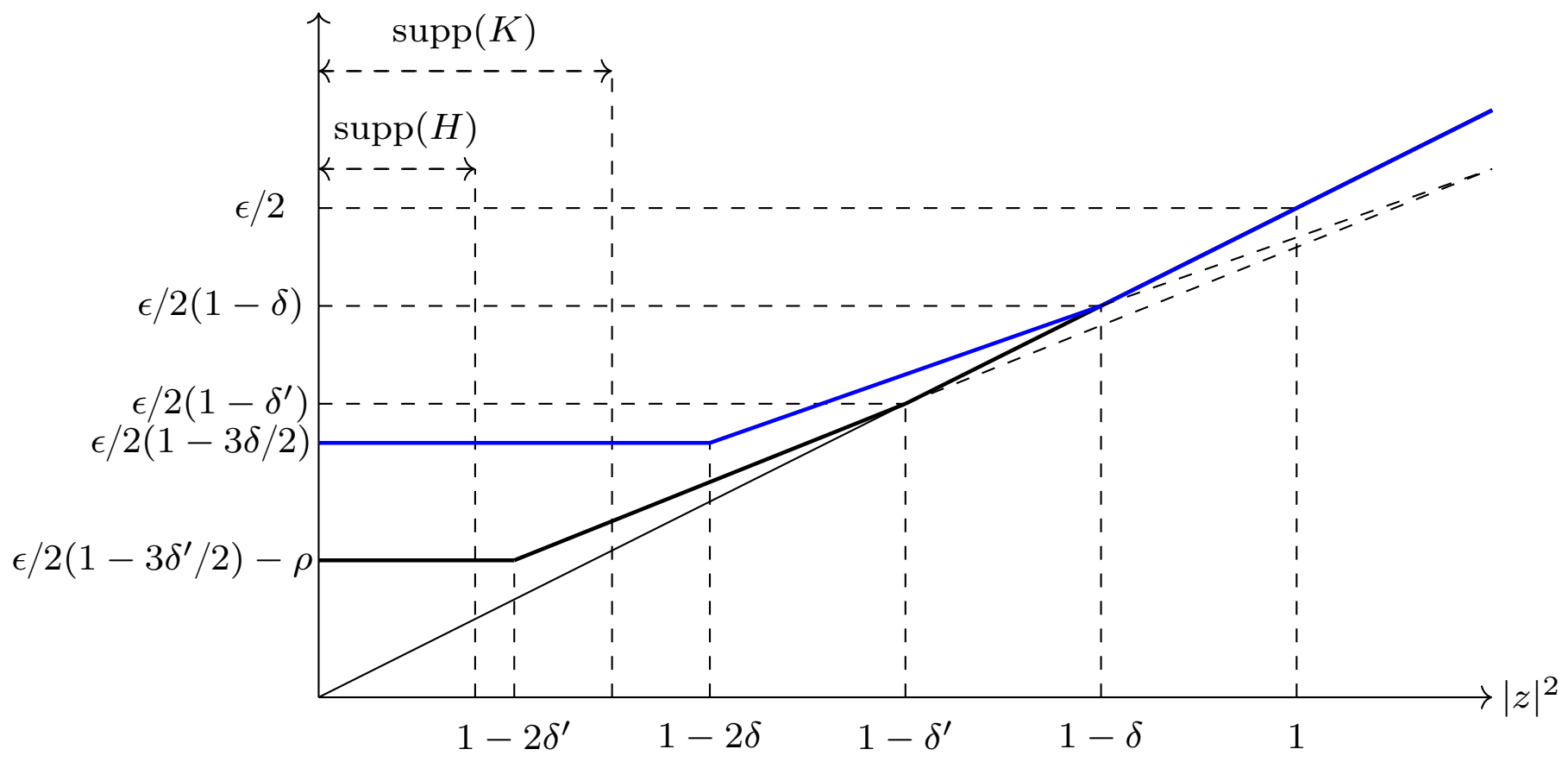
\begin{figure}[ht]
  \centering
  \begin{tikzpicture}[scale = 1.5]
\draw[->] (0,0) -- (0,3.5);
\draw[->] (0,0) -- (6,0);
\draw (0,0) -- (6,3);

\draw[dashed] (5,0) -- (5,2.5);
\draw[dashed] (0,2.5) -- (5,2.5);
\node at (5,-0.2) {\tiny{$1$}};

\draw[dashed] (4,0) -- (4,2);
\draw[dashed] (0,2) -- (4,2);
\node at (4,-0.2) {\tiny{$1 - \delta$}};

\draw[dashed] (3,0) -- (3,1.5);
\draw[dashed] (0,1.5) -- (3,1.5);
\node at (3,-0.2) {\tiny{$1 - \delta'$}};
\node at (-0.5,1.5) {\tiny{$\epsilon/2(1 - \delta')$}};

\draw[dashed] (2,0) -- (2,1.3);
\draw[dashed] (0,1.3) -- (2,1.3);
\node at (2,-0.2) {\tiny{$1 - 2\delta$}};

\draw[dashed] (1,0) -- (1,0.7);
\draw[dashed] (0,0.7) -- (1,0.7);
\node at (1,-0.2) {\tiny{$1 - 2\delta'$}};

\draw[dashed] (1.5,0) -- (1.5,3.2);
\draw[dashed, ->] (0,3.2) -- (1.5,3.2);
\draw[dashed, <-] (0,3.2) -- (1.5,3.2);
\node at (0.75,3.4) {\tiny{$\mathrm{supp}(K)$}};

\draw[dashed] (0.8,0) -- (0.8,2.7);
\draw[dashed, ->] (0,2.7) -- (0.8,2.7);
\draw[dashed, <-] (0,2.7) -- (0.8,2.7);
\node at (0.45,2.9) {\tiny{$\mathrm{supp}(H)$}};

\draw[line width=0.3mm, color = black] (3,1.5) -- (6,3);
\draw[line width=0.3mm, color = black] (1,0.7) -- (3,1.5);
\draw[line width=0.3mm, color = black] (0,0.7) -- (1,0.7);
\draw[dashed] (3,1.5) -- (6,2.7);
\node at (-0.8,0.7) {\tiny{$\epsilon/2(1 - 3\delta'/2) - \rho$}};

\draw[line width=0.3mm, color = blue] (4,2) -- (6,3);
\draw[line width=0.3mm, color = blue] (2,1.3) -- (4,2);
\draw[line width=0.3mm, color = blue] (0,1.3) -- (2,1.3);
\draw[dashed] (4,2) -- (6,2.7);
\node at (-0.6,1.3) {\tiny{$\epsilon/2(1 - 3\delta/2)$}};

\node at (6.2,0) {\tiny{$|z|^{2}$}};
\node at (-0.3, 2.5) {\tiny{$\epsilon/2$}};
\node at (-0.5, 2) {\tiny{$\epsilon/2(1-\delta)$}};

\end{tikzpicture}
\caption{A non-degenerate perturbation of $\phi \circ \phi^{1}_{\epsilon, \delta}$. The graphs of $H_{\epsilon/2, \delta}$ and $H'_{\epsilon/2, \delta'}$ are smoothening of the top and bottom curves near the "corners" respectively. Notice that $H'_{\epsilon/2,\delta'}$ is a slight modification of $H_{\epsilon/2, \delta'}$ by the term $-\rho$, compare with Figure \ref{fig:hamiltonian}; the small negative term $-\rho$ is there to create a slightly larger slope for $H'_{\epsilon/2,\delta'}$ compared to that of $H_{\epsilon/2,\delta'}$ in the region $\{1-2\delta' < |z|^{2} \leq 1-\delta'\}$.}
\label{fig:pert}
\end{figure}

Recall that $\phi$ is supported in $\{|z|^{2} < 1 - 2\delta'\}$ as it is generated by $H$. We now describe the perturbation $\psi_{\epsilon}$. Let us focus on the region $\{1- 2\delta \leq |z|^{2}\}$. On this region, 
$$\psi_{\epsilon} \text{ is generated by $H\# K\#$ (bottom graph }+\text{ top graph)}.$$
Define $\phi' = \psi_{\epsilon} \circ (\phi^{1}_{\epsilon, \delta})^{-1}$ which is supported in $\{|z|^{2} \leq 1 -\delta\}$ and,
$$\phi'\text{ is generated by $H \#K\#$ (bottom graph } - \text{ top graph)},$$
on the region $\{1-2\delta \leq |z|^{2} < 1 -\delta\}$. These must satisfy the following:
\begin{enumerate}
\item $\psi_{\epsilon}$ is rotation by angle $2\pi \epsilon$ on $\{1 - \delta \leq |z|^{2}\}$,
\item $\phi'$ is identity on $\{1 - \delta \leq |z|^{2}\}$,
\item $\psi_{\epsilon}$ is rotation by angles in $(0,2\pi \epsilon)$ on $\{1 - 2\delta \leq |z|^{2} < 1- \delta\}$,
\item $\phi'$ is rotation by angles in $(0,\pi \epsilon/2)$ on $\{1 - 2\delta \leq |z|^{2} < 1- \delta\}$. In particular, $\phi'$ has no orbits on this region.
\end{enumerate}
Let $\epsilon < \theta < 1$ be any other angle. Define, 
$$\psi_{\theta} := \psi_{\epsilon} \circ \phi^{1}_{\theta - \epsilon, \delta} = \phi' \circ \phi^{1}_{\theta, \delta},$$ 
which is generated by
\begin{equation}\label{eq:gen_ham_theta}
(\underbrace{H \# H_{\epsilon/2,\delta} \# H'_{\epsilon/2,\delta'} \# K}_\text{generator of $\psi_{\epsilon}$}) \# H_{\theta-\epsilon, \delta}.
\end{equation}
On the region $\{1-2\delta \leq |z|^{2}\},$ this Hamiltonian is equal to
$$\underbrace{\text{(bottom graph}-\text{top graph)}}_{\text{generator of } \phi'}+ H_{\theta,\delta}.$$
A quick calculation shows that:
\begin{enumerate}
\item $\psi_{\theta}$ is rotation by angle $2\pi \theta$ on $\{1 - \delta \leq |z|^{2}\}$.
\item $\psi_{\theta}$ is rotation by angles in $(0,2\pi \theta)$ on $\{1 - 2\delta \leq |z|^{2} < 1- \delta\}$.
\end{enumerate}
It is evident that $\psi_{\epsilon} = \psi_{\theta} = \phi'$ on the region $\{|z|^{2} < 1- 2\delta\},$ and hence
$$\mathrm{Fix}(\psi_{\epsilon}) = \mathrm{Fix}(\psi_{\theta}).$$
Thus, $\psi_{\theta}$ is a non-degenerate perturbation of $\phi \circ \phi^{1}_{\theta, \delta}$, the time-one map of $H \# H_{\theta, \delta}$, and has exactly the same fixed points as $\psi_{\epsilon}$. 

The filtered Floer complexes of $\psi_{\epsilon}$ and $\psi_{\theta}$ are identified because the two generating Hamiltonians \eqref{eq:gen_ham_eps} and \eqref{eq:gen_ham_theta} of $\psi_{\epsilon}$ and $\psi_{\theta}$ differ by a constant on $\{|z|^{2} \leq 1 - 2\delta\}$ where all the fixed points lie. To be more precise, one can apply the maximum principle \cite[Corollary 2.3]{Frauenfelder-Schlenk-07} to show that no Floer cylinder, neither for $\psi_{\epsilon}$ nor for $\psi_{\theta}$, can escape the region $\{1-2\delta \leq |z|^{2}\}$. And so not only the fixed points are the same but also the Floer differentials are identified. However, the actions of fixed points are all shifted by a constant.

Since the boundary depth is invariant under constant shifts of actions,
$$\beta(\psi_{\theta}, \mathbb{K}) = \beta(\psi_{\epsilon}, \mathbb{K}).$$ 
The maps $\psi_{\epsilon}$ and $\psi_{\theta}$ only depend on the perturbation datum $K, \delta'$. Let us denote the boundary depth of $\psi_{\epsilon}$ by $\beta(\phi, \epsilon, K, \delta')$ and the boundary depth of $\psi_{\theta}$ by $\beta(\phi,\theta, K, \delta')$. Then the following holds:
$$\beta(\phi \circ \phi^{1}_{\epsilon, \delta}) = \lim_{\substack{K \to 0\\ \delta' \to \delta}} \beta(\phi,\epsilon, K, \delta') = \lim_{\substack{K \to 0\\ \delta' \to \delta}} \beta(\phi,\theta, K, \delta') = \beta(\phi \circ \phi^{1}_{\theta, \delta}),$$
This shows the independence of the angle at infinity. Independence of $\delta$ is implicit in the proof above. Consider the two Hamiltonians $H\# H_{\eps,\delta}$ and $H \# H_{\eps, \delta'}$ for some $1>\eps >0$ and small numbers $\delta'> \delta > 0$. Take some number $\delta'' > \delta'$ and close enough to $\delta'$. Then consider the non-degenerate perturbations of the following form defined similarly as in Figure \ref{fig:pert}:
\begin{enumerate}
\item $H \# H_{\eps/2,\delta}\#H'_{\eps/2,\delta''}\# K$ for $H \# H_{\eps, \delta},$
\item $H \# H_{\eps/2,\delta'}\#H'_{\eps/2,\delta''}\# K$ for $H \# H_{\eps, \delta'},$
\end{enumerate}
where $K$ is a $C^{2}$-small compactly supported Hamiltonian in $|z| < 1-2\delta'$. Then, similar analysis as above shows that the two perturbations have no fixed points on $|z|> 1-2\delta''$, moreover inside $\{|z|<1-2\delta''\}$, where the fixed points lie, they differ by a constant, which does not affect the boundary depth.
\end{proof}

\bibliographystyle{amsalpha}
\bibliography{merging}

@article{franks2013entropy,
  title={Entropy zero area preserving diffeomorphisms of $S^2$},
  author={Franks, John and Handel, Michael},
  journal={Geometry \& Topology},
  volume={16},
  pages={2187--2284},
  year={2013},
  publisher={Mathematical Sciences Publishers}
}

@article{le2022topological,
  title={Topological horseshoes for surface homeomorphisms},
  author={Le Calvez, Patrice and Tal, Fabio},
  journal={Duke Mathematical Journal},
  volume={171},
  number={12},
  pages={2519--2626},
  year={2022},
  publisher={Duke University Press}
}

@article{bramham2024spectral,
  title={Spectral invariants for non-compactly supported Hamiltonians on the disc, and an application to the mean action spectrum},
  author={Bramham, Barney and Pirnapasov, Abror},
  journal={arXiv preprint arXiv:2403.07863},
  year={2024},
  pages={1--15},
}

@article{viterbo1992symplectic,
  title={Symplectic topology as the geometry of generating functions},
  author={Viterbo, Claude},
  journal={Mathematische Annalen},
  volume={292},
  number={1},
  pages={685--710},
  year={1992}
}

@phdthesis{senior2023asymptotic,
  title={The Asymptotic Action and Generalizations of a Theorem of Fathi},
  author={Senior, David Bechara},
  year={2023},
  school={Ruhr-Universit{\"a}t Bochum}
}

@article{whitney1934differentiable,
  title={Differentiable functions defined in closed sets. I},
  author={Whitney, Hassler},
  journal={Transactions of the American Mathematical Society},
  volume={36},
  number={2},
  pages={369--387},
  year={1934},
  publisher={JSTOR}
}

@article{hutchings2025reeb,
  title={Reeb orbits frequently intersecting a symplectic surface},
  author={Hutchings, Michael},
  journal={arXiv preprint arXiv:2504.19332},
  year={2025},
pages={1--37}
}

@article{abbondandolo2018sharp,
  title={Sharp systolic inequalities for Reeb flows on the three-sphere},
  author={Abbondandolo, Alberto and Bramham, Barney and Hryniewicz, Umberto L and Salom{\~a}o, Pedro AS},
  journal={Inventiones mathematicae},
  volume={211},
  number={2},
  pages={687--778},
  year={2018},
  publisher={Springer}
}

@article{pirnapasov2024hutchings,
  title={Hutchings’ inequality for the {C}alabi invariant revisited with an application to pseudo-rotations},
  author={Pirnapasov, Abror},
  journal={Journal of Fixed Point Theory and Applications},
  volume={26},
  number={4},
  pages={56},
  year={2024},
  publisher={Springer}
}

@article{herman1979conjugaison,
  title={Sur la conjugaison diff{\'e}rentiable des diff{\'e}omorphismes du cercle {\`a} des rotations},
  author={Herman, Michael R},
  journal={Publications Math{\'e}matiques de l'IH{\'E}S},
  volume={49},
  pages={5--233},
  year={1979}
}

@inproceedings{yoccoz1984conjugaison,
  title={Conjugaison diff{\'e}rentiable des diff{\'e}omorphismes du cercle dont le nombre de rotation v{\'e}rifie une condition diophantienne},
  author={Yoccoz, J-C},
  booktitle={Annales scientifiques de l'{\'E}cole Normale Sup{\'e}rieure},
  volume={17},
  number={3},
  pages={333--359},
  year={1984}
}

@article{banyaga1974formes,
  title={Formes-volume sur les vari{\'e}t{\'e}sa bord},
  author={Banyaga, Augustin},
  journal={Enseignement Math},
  volume={20},
  number={2},
  pages={127--131},
  year={1974}
}

@article{franks1988generalizations,
  title={Generalizations of the Poincar{\'e}-Birkhoff theorem},
  author={Franks, John},
  journal={Annals of Mathematics},
  volume={128},
  number={1},
  pages={139--151},
  year={1988},
  publisher={JSTOR}
}

@article{alves2021braid,
  title={Braid stability and the Hofer metric},
  author={Alves, Marcelo RR and Meiwes, Matthias},
  journal={arXiv preprint arXiv:2112.11351},
  year={2021},
  PAGES = {1--56},
}

@article{kislev2022bounds,
  title={Bounds on spectral norms and barcodes},
  author={Kislev, Asaf and Shelukhin, Egor},
  journal={Geometry \& Topology},
  volume={25},
  number={7},
  pages={3257--3350},
  year={2022},
  publisher={Mathematical Sciences Publishers}
}

@article{shelukhin-zhao-22,
  title={The $\backslash$mathbb $\mathbb {Z}/p\mathbb {Z}$-equivariant product-isomorphism in fixed point Floer cohomology},
  author={Shelukhin, Egor and Zhao, Jingyu},
  journal={Journal of Symplectic Geometry},
  volume={19},
  number={5},
  pages={1101--1188},
  year={2022},
  publisher={International Press of Boston}
}

@article{katok1980lyapunov,
  title={Lyapunov exponents, entropy and periodic orbits for diffeomorphisms},
  author={Katok, Anatole},
  journal={Publications Math{\'e}matiques de l'IH{\'E}S},
  volume={51},
  pages={137--173},
  year={1980}
}

@article{bechara2023asymptotic,
  title={The asymptotic mean action and the asymptotic linking number for pseudo-rotations},
  author={Senior, David Bechara and Calvez, Patrice Le and Pirnapasov, Abror},
  journal={arXiv preprint arXiv:2312.17087},
  year={2023},
pages ={1--20}
}

@article {floer-hofer-salamon-transversality-95,
    AUTHOR = {Floer, A. and Hofer, H. and Salamon, D.},
     TITLE = {Transversality in elliptic {M}orse theory for the symplectic
              action},
   JOURNAL = {Duke Math. J.},
  FJOURNAL = {Duke Mathematical Journal},
    VOLUME = {80},
      YEAR = {1995},
    NUMBER = {1},
     PAGES = {251--292},
      ISSN = {0012-7094},
   MRCLASS = {58E05 (57R70 58F05)},
  MRNUMBER = {1360618},
MRREVIEWER = {Karl Friedrich Siburg},
       DOI = {10.1215/S0012-7094-95-08010-7},
       URL = {https://doi.org/10.1215/S0012-7094-95-08010-7},
}

@article {Viterbo-99,
    AUTHOR = {Viterbo, C.},
     TITLE = {Functors and computations in {F}loer homology with
              applications. {I}},
   JOURNAL = {Geom. Funct. Anal.},
  FJOURNAL = {Geometric and Functional Analysis},
    VOLUME = {9},
      YEAR = {1999},
    NUMBER = {5},
     PAGES = {985--1033},
      ISSN = {1016-443X},
   MRCLASS = {53D40 (57R17 57R58)},
  MRNUMBER = {1726235},
MRREVIEWER = {David E. Hurtubise},
       DOI = {10.1007/s000390050106},
       URL = {https://doi.org/10.1007/s000390050106},
}

@incollection {ginzburg-05,
    AUTHOR = {Ginzburg, V. L.},
     TITLE = {The {W}einstein conjecture and theorems of nearby and almost
              existence},
 BOOKTITLE = {The breadth of symplectic and {P}oisson geometry},
    VOLUME = {232},
     PAGES = {139--172},
 PUBLISHER = {Birkh\"{a}user Boston, Boston, MA},
      YEAR = {2005},
   MRCLASS = {53D40 (37J05 37J45)},
  MRNUMBER = {2103006},
MRREVIEWER = {Kai Cieliebak},
       DOI = {10.1007/0-8176-4419-9\_6},
       URL = {https://doi.org/10.1007/0-8176-4419-9_6},
}

@article {lanzat-13,
    AUTHOR = {Lanzat, S.},
     TITLE = {Quasi-morphisms and symplectic quasi-states for convex
              symplectic manifolds},
   JOURNAL = {Int. Math. Res. Not. IMRN},
  FJOURNAL = {International Mathematics Research Notices. IMRN},
      YEAR = {2013},
    NUMBER = {23},
     PAGES = {5321--5365},
      ISSN = {1073-7928},
   MRCLASS = {53D40},
  MRNUMBER = {3142258},
MRREVIEWER = {Michael J. Usher},
       DOI = {10.1093/imrn/rns205},
       URL = {https://doi.org/10.1093/imrn/rns205},
}

@article {Lanzat-16,
    AUTHOR = {Lanzat, S.},
     TITLE = {Hamiltonian {F}loer homology for compact convex symplectic
              manifolds},
   JOURNAL = {Beitr. Algebra Geom.},
  FJOURNAL = {Beitr\"{a}ge zur Algebra und Geometrie. Contributions to Algebra
              and Geometry},
    VOLUME = {57},
      YEAR = {2016},
    NUMBER = {2},
     PAGES = {361--390},
      ISSN = {0138-4821},
   MRCLASS = {53D05 (53D40 53D45)},
  MRNUMBER = {3493995},
MRREVIEWER = {Umberto Leone Hryniewicz},
       DOI = {10.1007/s13366-015-0254-6},
       URL = {https://doi.org/10.1007/s13366-015-0254-6},
}

@article {alizadeh-atallah-cant,
    AUTHOR = {Alizadeh, H. and Atallah, M. S. and Cant, D.},
     TITLE = {Lagrangian intersections and the spectral norm in
              convex-at-infinity symplectic manifolds},
   JOURNAL = {Math. Z.},
  FJOURNAL = {Mathematische Zeitschrift},
    VOLUME = {310},
      YEAR = {2025},
    NUMBER = {1},
     PAGES = {2--48},
      ISSN = {0025-5874},
   MRCLASS = {53D12 (37J06 53D40)},
  MRNUMBER = {4879080},
       DOI = {10.1007/s00209-025-03710-0},
       URL = {https://doi.org/10.1007/s00209-025-03710-0},
}

@article{alizadeh-atallah-cant-ellipsoid-24,
  title={The spectral diameter of a symplectic ellipsoid},
  author={Alizadeh, H. and Atallah, M. S and Cant, D.},
  journal={Journal of Symplectic Geometry},
  year={2024},
pages={1--21}
}

@article {Usher-13,
    AUTHOR = {Usher, M.},
     TITLE = {Hofer's metrics and boundary depth},
   JOURNAL = {Ann. Sci. \'{E}c. Norm. Sup\'{e}r. (4)},
  FJOURNAL = {Annales Scientifiques de l'\'{E}cole Normale Sup\'{e}rieure. Quatri\`eme
              S\'{e}rie},
    VOLUME = {46},
      YEAR = {2013},
     PAGES = {57--128},
      ISSN = {0012-9593},
   MRCLASS = {53D40 (53D12 53D35)},
  MRNUMBER = {3087390},
MRREVIEWER = {David E. Hurtubise},
       DOI = {10.24033/asens.2185},
       URL = {https://doi.org/10.24033/asens.2185},
}

@article{feng-zhang,
  title={Spectrally-large scale geometry in cotangent bundles},
  author={Feng, Q. and Zhang, J.},
  journal={arXiv preprint arXiv:2401.17590},
  year={2024},
pages = {1--35}
}

@article {Schwarz-00,
    AUTHOR = {Schwarz, M.},
     TITLE = {On the action spectrum for closed symplectically aspherical
              manifolds},
   JOURNAL = {Pacific J. Math.},
  FJOURNAL = {Pacific Journal of Mathematics},
    VOLUME = {193},
      YEAR = {2000},
    NUMBER = {2},
     PAGES = {419--461},
      ISSN = {0030-8730},
   MRCLASS = {53D40 (57R58)},
  MRNUMBER = {1755825},
MRREVIEWER = {David E. Hurtubise},
       DOI = {10.2140/pjm.2000.193.419},
       URL = {https://doi.org/10.2140/pjm.2000.193.419},
}

@article {Seyfaddini-15,
    AUTHOR = {Seyfaddini, S.},
     TITLE = {Spectral killers and {P}oisson bracket invariants},
   JOURNAL = {J. Mod. Dyn.},
  FJOURNAL = {Journal of Modern Dynamics},
    VOLUME = {9},
      YEAR = {2015},
     PAGES = {51--66},
      ISSN = {1930-5311},
   MRCLASS = {53D40 (37J05)},
  MRNUMBER = {3395260},
MRREVIEWER = {Chun-Gen Liu},
       DOI = {10.3934/jmd.2015.9.51},
       URL = {https://doi.org/10.3934/jmd.2015.9.51},
}

@article {Frauenfelder-Schlenk-07,
    AUTHOR = {Frauenfelder, U. and Schlenk, F.},
     TITLE = {Hamiltonian dynamics on convex symplectic manifolds},
   JOURNAL = {Israel J. Math.},
  FJOURNAL = {Israel Journal of Mathematics},
    VOLUME = {159},
      YEAR = {2007},
     PAGES = {1--56},
      ISSN = {0021-2172},
   MRCLASS = {53D40 (37J05 37J45)},
  MRNUMBER = {2342472},
MRREVIEWER = {Michael J. Usher},
       DOI = {10.1007/s11856-007-0037-3},
       URL = {https://doi.org/10.1007/s11856-007-0037-3},
}

@incollection {oh-spec-inv-05,
    AUTHOR = {Oh, Yong-Geun},
     TITLE = {Construction of spectral invariants of {H}amiltonian paths on
              closed symplectic manifolds},
 BOOKTITLE = {The breadth of symplectic and {P}oisson geometry},
    VOLUME = {232},
     PAGES = {525--570},
 PUBLISHER = {Birkh\"{a}user Boston, Boston, MA},
      YEAR = {2005},
      ISBN = {0-8176-3565-3},
   MRCLASS = {53D40 (53D45 58E30)},
  MRNUMBER = {2103018},
MRREVIEWER = {Darko\ Milinkovi\'{c}},
       DOI = {10.1007/0-8176-4419-9\{_}18}

@article {ginzburg-gurel-hyperbolic,
    AUTHOR = {Ginzburg, V. L. and G\"{u}rel, B. Z.},
     TITLE = {Hyperbolic fixed points and periodic orbits of {H}amiltonian
              diffeomorphisms},
   JOURNAL = {Duke Math. J.},
  FJOURNAL = {Duke Mathematical Journal},
    VOLUME = {163},
      YEAR = {2014},
    NUMBER = {3},
     PAGES = {565--590},
      ISSN = {0012-7094,1547-7398},
   MRCLASS = {53D40 (37J10 70H12)},
  MRNUMBER = {3165423},
MRREVIEWER = {Cristian\ L\u{a}zureanu},
       DOI = {10.1215/00127094-2410433},
       URL = {https://doi.org/10.1215/00127094-2410433},
}

@article {ginzburg-gurel-invent,
    AUTHOR = {Ginzburg, V. L. and G\"{u}rel, B. Z.},
     TITLE = {Hamiltonian pseudo-rotations of projective spaces},
   JOURNAL = {Invent. Math.},
  FJOURNAL = {Inventiones Mathematicae},
    VOLUME = {214},
      YEAR = {2018},
    NUMBER = {3},
     PAGES = {1081--1130},
      ISSN = {0020-9910,1432-1297},
   MRCLASS = {53D40 (37J10 37J45)},
  MRNUMBER = {3878727},
       DOI = {10.1007/s00222-018-0818-9},
       URL = {https://doi.org/10.1007/s00222-018-0818-9},
}

@article {milinkovic-oh-97,
    AUTHOR = {Milinkovi\'{c}, D. and Oh, Y-G},
     TITLE = {Floer homology as the stable {M}orse homology},
   JOURNAL = {J. Korean Math. Soc.},
  FJOURNAL = {Journal of the Korean Mathematical Society},
    VOLUME = {34},
      YEAR = {1997},
    NUMBER = {4},
     PAGES = {1065--1087},
      ISSN = {0304-9914,2234-3008},
   MRCLASS = {58E05 (58F05)},
  MRNUMBER = {1485969},
MRREVIEWER = {Joa\ Weber},
}

@article {Handel-90-rot-set-is-closed-int,
    AUTHOR = {Handel, M.},
     TITLE = {The rotation set of a homeomorphism of the annulus is closed},
   JOURNAL = {Comm. Math. Phys.},
  FJOURNAL = {Communications in Mathematical Physics},
    VOLUME = {127},
      YEAR = {1990},
    NUMBER = {2},
     PAGES = {339--349},
      ISSN = {0010-3616,1432-0916},
   MRCLASS = {58F11 (58F18)},
  MRNUMBER = {1037109},
MRREVIEWER = {Edward\ Slaminka},
       URL = {http://projecteuclid.org/euclid.cmp/1104180142},
}

@article {Bechara-23-ann-du-quebec,
    AUTHOR = {Bechara Senior, D.},
     TITLE = {Asymptotic action and asymptotic winding number for
              area-preserving diffeomorphisms of the disk},
   JOURNAL = {Ann. Math. Qu\'{e}.},
  FJOURNAL = {Annales Math\'{e}matiques du Qu\'{e}bec},
    VOLUME = {47},
      YEAR = {2023},
    NUMBER = {2},
     PAGES = {443--459},
      ISSN = {2195-4755,2195-4763},
   MRCLASS = {37E30 (37C05)},
  MRNUMBER = {4645699},
MRREVIEWER = {Jonathan\ Conejeros},
       DOI = {10.1007/s40316-021-00173-7},
       URL = {https://doi.org/10.1007/s40316-021-00173-7},
}

@article {salamon-zehnder-92,
    AUTHOR = {Salamon, D. and Zehnder, E.},
     TITLE = {Morse theory for periodic solutions of {H}amiltonian systems
              and the {M}aslov index},
   JOURNAL = {Comm. Pure Appl. Math.},
  FJOURNAL = {Communications on Pure and Applied Mathematics},
    VOLUME = {45},
      YEAR = {1992},
    NUMBER = {10},
     PAGES = {1303--1360},
      ISSN = {0010-3640,1097-0312},
   MRCLASS = {58E05 (34C25 58F22 70H05)},
  MRNUMBER = {1181727},
MRREVIEWER = {Yong-Geun\ Oh},
       DOI = {10.1002/cpa.3160451004},
       URL = {https://doi.org/10.1002/cpa.3160451004},
}

@article {fish-11-gromov-comp,
    AUTHOR = {Fish, J. W.},
     TITLE = {Target-local {G}romov compactness},
   JOURNAL = {Geom. Topol.},
  FJOURNAL = {Geometry \& Topology},
    VOLUME = {15},
      YEAR = {2011},
    NUMBER = {2},
     PAGES = {765--826},
      ISSN = {1465-3060,1364-0380},
   MRCLASS = {53D45 (32Q65 53C15)},
  MRNUMBER = {2800366},
MRREVIEWER = {Antoine\ F. R. Gournay},
       DOI = {10.2140/gt.2011.15.765},
       URL = {https://doi.org/10.2140/gt.2011.15.765},
}

@book {siburg-princ-least-act-04,
    AUTHOR = {Siburg, K. F.},
     TITLE = {The principle of least action in geometry and dynamics},
    VOLUME = {1844},
 PUBLISHER = {Springer-Verlag, Berlin},
      YEAR = {2004},
     PAGES = {xii+128},
      ISBN = {3-540-21944-7},
   MRCLASS = {37J50 (37E40 37J10 53D35 58E30)},
  MRNUMBER = {2076302},
MRREVIEWER = {Renato\ Iturriaga},
       DOI = {10.1007/b97327},
       URL = {https://doi.org/10.1007/b97327},
}

@article {MR1581099,
    AUTHOR = {Baer, R.},
     TITLE = {Kurventypen auf {F}l\"{a}chen},
   JOURNAL = {J. Reine Angew. Math.},
  FJOURNAL = {Journal f\"{u}r die Reine und Angewandte Mathematik. [Crelle's
              Journal]},
    VOLUME = {156},
      YEAR = {1927},
     PAGES = {231--246},
      ISSN = {0075-4102,1435-5345},
   MRCLASS = {DML},
  MRNUMBER = {1581099},
       DOI = {10.1515/crll.1927.156.231},
       URL = {https://doi.org/10.1515/crll.1927.156.231},
}

@article {MR1581156,
    AUTHOR = {Baer, R.},
     TITLE = {Isotopie von {K}urven auf orientierbaren, geschlossenen
              {F}l\"{a}chen und ihr {Z}usammenhang mit der topologischen
              {D}eformation der {F}l\"{a}chen},
   JOURNAL = {J. Reine Angew. Math.},
  FJOURNAL = {Journal f\"{u}r die Reine und Angewandte Mathematik. [Crelle's
              Journal]},
    VOLUME = {159},
      YEAR = {1928},
     PAGES = {101--116},
      ISSN = {0075-4102,1435-5345},
   MRCLASS = {DML},
  MRNUMBER = {1581156},
       DOI = {10.1515/crll.1928.159.101},
       URL = {https://doi.org/10.1515/crll.1928.159.101},
}

@article {MR214087,
    AUTHOR = {Epstein, D. B. A.},
     TITLE = {Curves on {$2$}-manifolds and isotopies},
   JOURNAL = {Acta Math.},
  FJOURNAL = {Acta Mathematica},
    VOLUME = {115},
      YEAR = {1966},
     PAGES = {83--107},
      ISSN = {0001-5962,1871-2509},
   MRCLASS = {57.20},
  MRNUMBER = {214087},
MRREVIEWER = {P.\ Dedecker},
       DOI = {10.1007/BF02392203},
       URL = {https://doi.org/10.1007/BF02392203},
}

@book {MR2850125,
    AUTHOR = {Farb, Benson and Margalit, Dan},
     TITLE = {A primer on mapping class groups},
    SERIES = {Princeton Mathematical Series},
    VOLUME = {49},
 PUBLISHER = {Princeton University Press, Princeton, NJ},
      YEAR = {2012},
     PAGES = {xiv+472},
      ISBN = {978-0-691-14794-9},
   MRCLASS = {57M50 (20F36 20F65 57M07 57N05)},
  MRNUMBER = {2850125},
MRREVIEWER = {Stephen\ P.\ Humphries},
}

@article {LeC-duke,
    AUTHOR = {Le Calvez, P.},
     TITLE = {Periodic orbits of {H}amiltonian homeomorphisms of surfaces},
   JOURNAL = {Duke Math. J.},
  FJOURNAL = {Duke Mathematical Journal},
    VOLUME = {133},
      YEAR = {2006},
    NUMBER = {1},
     PAGES = {125--184},
      ISSN = {0012-7094,1547-7398},
   MRCLASS = {37E30 (37C25 37J45)},
  MRNUMBER = {2219272},
MRREVIEWER = {Shigenori\ Matsumoto},
       DOI = {10.1215/S0012-7094-06-13315-X},
       URL = {https://doi.org/10.1215/S0012-7094-06-13315-X},
}

@article {P-02,
    AUTHOR = {Polterovich, L.},
     TITLE = {Growth of maps, distortion in groups and symplectic geometry},
   JOURNAL = {Invent. Math.},
  FJOURNAL = {Inventiones Mathematicae},
    VOLUME = {150},
      YEAR = {2002},
    NUMBER = {3},
     PAGES = {655--686},
      ISSN = {0020-9910,1432-1297},
   MRCLASS = {53D35 (53D40)},
  MRNUMBER = {1946555},
MRREVIEWER = {Karl\ Friedrich\ Siburg},
       DOI = {10.1007/s00222-002-0251-x},
       URL = {https://doi.org/10.1007/s00222-002-0251-x},
}

@article{Schwarz2000,
  author  = {Schwarz, Matthias},
  title   = {On the action spectrum for closed symplectically
             aspherical manifolds},
  journal = {Pacific Journal of Mathematics},
  volume  = {193},
  number  = {2},
  pages   = {419--461},
  year    = {2000},
  doi     = {10.2140/pjm.2000.193.419}
}

@incollection{Oh2005,
  author    = {Oh, Yong-Geun},
  title     = {Construction of spectral invariants of Hamiltonian
               paths on closed symplectic manifolds},
  booktitle = {The Breadth of Symplectic and Poisson Geometry},
  series    = {Progress in Mathematics},
  volume    = {232},
  pages     = {525--570},
  publisher = {Birkh{\"a}user},
  address   = {Boston},
  year      = {2005}
}

@article {MR1511684,
    AUTHOR = {Brouwer, L. E. J.},
     TITLE = {Beweis des ebenen {T}ranslationssatzes},
   JOURNAL = {Math. Ann.},
  FJOURNAL = {Mathematische Annalen},
    VOLUME = {72},
      YEAR = {1912},
    NUMBER = {1},
     PAGES = {37--54},
      ISSN = {0025-5831,1432-1807},
   MRCLASS = {DML},
  MRNUMBER = {1511684},
       DOI = {10.1007/BF01456888},
       URL = {https://doi.org/10.1007/BF01456888},
}

@article {MR925994,
    AUTHOR = {Fathi, Albert},
     TITLE = {An orbit closing proof of {B}rouwer's lemma on translation
              arcs},
   JOURNAL = {Enseign. Math. (2)},
  FJOURNAL = {L'Enseignement Math\'{e}matique. Revue Internationale. 2e
              S\'{e}rie},
    VOLUME = {33},
      YEAR = {1987},
    NUMBER = {3-4},
     PAGES = {315--322},
      ISSN = {0013-8584},
   MRCLASS = {55M20 (54H20 54H25 58F25)},
  MRNUMBER = {925994},
MRREVIEWER = {Edward\ Slaminka},
}

@book {MR2193309,
    AUTHOR = {Milnor, J.},
     TITLE = {Dynamics in one complex variable},
    SERIES = {Annals of Mathematics Studies},
    VOLUME = {160},
   EDITION = {Third},
 PUBLISHER = {Princeton University Press, Princeton, NJ},
      YEAR = {2006},
     PAGES = {viii+304},
      ISBN = {978-0-691-12488-9; 0-691-12488-4},
   MRCLASS = {37Fxx (30-01 30D05 37-01)},
  MRNUMBER = {2193309},
}
\end{document}